\documentclass[11pt,letterpaper]{article}
\usepackage[T1]{fontenc}
\usepackage[utf8]{inputenc}
\usepackage[margin=1in]{geometry}
\usepackage{times}
\usepackage{mathtools,amssymb,amsthm}
\usepackage{booktabs,microtype,array}
\usepackage[authoryear,round]{natbib}
\usepackage{xurl}
\usepackage{mleftright}
\usepackage[hidelinks]{hyperref}
\mleftright
\hypersetup{pdftitle={Sharp Critical Rates for Locally Private Distribution Estimation},pdfauthor={Yi Liu, Xiaotian Chang, Qirui Hu}}
\newcommand{\E}{\mathbb E}
\newcommand{\Prob}{\mathbb P}
\newcommand{\R}{\mathbb R}
\newcommand{\T}{\mathbb T}
\DeclareMathOperator{\tr}{tr}
\DeclareMathOperator{\Var}{Var}
\DeclareMathOperator{\sign}{sign}
\numberwithin{equation}{section}
\newtheorem{theorem}{Theorem}
\newtheorem{lemma}[theorem]{Lemma}
\newtheorem{proposition}[theorem]{Proposition}
\newtheorem{corollary}[theorem]{Corollary}
\newenvironment{keywords}{\par\smallskip\noindent\textbf{Keywords:}\ }{\par\smallskip}
\bibpunct{(}{)}{;}{a}{,}{,}
\title{Sharp Critical Rates for Locally Private Distribution Estimation}
\author{Yi Liu\textsuperscript{1}\quad Xiaotian Chang\textsuperscript{2}\quad Qirui Hu\textsuperscript{3}\\[0.5em]
\small\textsuperscript{1}York University\\
\small\textsuperscript{2}Nanyang Technological University\\
\small\textsuperscript{3}Shanghai University of Finance and Economics}
\date{}

\begin{document}
\maketitle

\begin{abstract}
We determine sharp minimax rates for distribution estimation under pure local differential privacy at critical query smoothness. For the integral probability metric generated by the unit ball of $B^d_{2,q}$, $1\le q\le\infty$, over fixed periodic Sobolev density balls with positive smoothness and interior slack at the uniform density, the expected risk is of order $N^{-1/2}(\log N)^{3/2-1/q}$, where $N=n\varepsilon^2$; expected squared risk has the squared order. The fine index governs aggregation across resolutions that share one private information budget. A single noninteractive release with a measurable proper density output attains these rates, and matching lower bounds allow finite public blackboard interaction and repeated participation under each user's total privacy budget. A known $L^2$ radius permits the same rates without using density smoothness, with output in the $L^2$ density ball. Separately, for positive smooth densities on $[0,1]$, integrated squared CDF risk has minimax order $(\log N)^2/N$, with sharp quantile, Wasserstein, and transport map consequences. The rates hold uniformly for $0<\varepsilon\le1$ as $N\to\infty$, with class parameters fixed.
\end{abstract}


\section{Introduction}
\label{sec:introduction}

Under local differential privacy, each of $n$ users releases privatized information about one observation under a total privacy budget $\varepsilon$. The same reports support many expectation queries, whose accuracy depends on query smoothness, motivating our study of how estimation errors accumulate across spatial resolutions when both densities and queries are smooth. At the critical point where Sobolev query smoothness equals the ambient dimension $d$, the bounds of \citet{AlbertEtAl2026} leave a gap between $N^{-1/2}$ and $N^{-1/2}\left(\log N\right)^{2d}$, where $N=n\varepsilon^2$. We close this gap by establishing the sharp minimax expected risk $\left(\log N\right)/\sqrt N$, thereby determining the exact logarithmic factor at criticality. One noninteractive release attains this rate even against finite public blackboard interaction with repeated visits under each user's total budget. Varying the Besov fine index then reveals how aggregation across resolutions determines the logarithmic exponent.

More generally, for the integral probability metric (IPM) generated by the unit ball of $B^d_{2,q}$, we establish
\begin{equation*}
\inf_{Q,\widehat P}\sup_{f\in\mathcal F_\beta\left(R\right)}
\E_f d_q\left(\widehat P,P_f\right)
\asymp \frac{\left(\log N\right)^{3/2-1/q}}{\sqrt N},
\qquad 1\le q\le\infty,
\end{equation*}
where $d_q$ denotes this metric, $\mathcal F_\beta\left(R\right)$ is a periodic Sobolev density ball with fixed smoothness $\beta>0$ and radius $R>1$, and $P_f$ is the law with density $f$. The infimum ranges over the permitted public protocols $Q$ and arbitrary probability-measure estimates $\widehat P$. One noninteractive density output attains the rate for every $q$; expected squared loss has its squared order.

The logarithmic penalty at criticality reflects a multiscale information constraint: many relevant resolutions contribute at the same order, but must share a single private information budget. Its exponent is determined jointly by this constraint and how the query loss aggregates errors across resolutions. To capture this simultaneous difficulty, we choose dimension-matched prior widths that keep the entire lower-bound family inside a fixed smooth density class, and combine Fisher information contraction \citep{BarnesChenOzgur2020}, coordinatewise van Trees inequalities \citep{GillLevit1995}, and clipping to obtain the sharp first-moment lower bound. For the upper bound, ellipsoid scaling \citep{NikolovTang2024} and spherical randomization \citep{DuchiJordanWainwright2018} balance the weighted block errors, while stronger concentration in higher-dimensional blocks controls their maximum without an additional $\sqrt{\log\log N}$ factor.

For cumulative distribution functions, the queries are threshold probabilities, and we establish the minimax integrated squared CDF risk $\left(\log N\right)^2/N$ over fixed positive H\"older classes on $\left[0,1\right]$, matching the order attained by hierarchical range-query procedures \citep{CormodeKulkarniSrivastava2019}. Monotone projection tuned using only $n,\varepsilon$ transfers this rate to quantile and squared Wasserstein risks, while inversion with specified density bounds yields the sharp transport map rate. The latter resolves the one-dimensional local privacy gap of \citet{LalanneEtAl2026} on the class specified in Section~\ref{sec:cdf}.

Table~\ref{tab:closest-rates} summarizes the two rate comparisons.

\begin{table}[ht]
\centering
\small
\begin{tabular}{@{}>{\raggedright\arraybackslash}p{0.30\linewidth}ccc@{}}
\toprule
Problem and reference & Prior lower & Prior upper & \textbf{Minimax rate (this work)} \\
\midrule
Expected Sobolev IPM loss\newline \citet{AlbertEtAl2026}
& $N^{-1/2}$ & $N^{-1/2}\left(\log N\right)^{2d}$ & $N^{-1/2}\log N$ \\
\addlinespace
Squared transport map loss\newline \citet{LalanneEtAl2026}
& $N^{-1}$ & $\max\left\{n^{-1},N^{-\frac{s+1}{s+3/2}}\right\}$ & $N^{-1}\left(\log N\right)^2$ \\
\bottomrule
\end{tabular}
\caption{Rates for expected losses, with $N=n\varepsilon^2$ and fixed class parameters; $s$ denotes density smoothness. See Sections~\ref{sec:theory} and~\ref{sec:cdf} for the comparison conditions.}
\label{tab:closest-rates}
\end{table}

Besov IPMs were studied for nonprivate distribution estimation by \citet{UppalSinghPoczos2019}. In the locally private setting, our analysis identifies how a shared information budget interacts with the fine index of the query loss. The distinction between query losses is central to the comparison with existing private estimators. Strong-loss private density estimation \citep{ButuceaEtAl2020} treats a different loss geometry.

Discrete query-release lower bounds allow biased estimators and sequential interaction \citep{BlasiokEtAl2019}, while spectral workload bounds \citep{McKennaEtAl2020} concern linear unbiased decoding. Our priors connect the information budget to continuous losses within fixed smooth classes. Related information inequalities appear in interactive federated estimation \citep{CaiLi2026}. Factorization-norm bounds for noninteractive local linear queries \citep{EdmondsNikolovUllman2020} include finite-domain thresholds. The geometric framework of \citet{NikolovTang2024} treats unbiased mean estimation under central approximate and concentrated privacy, and also under noninteractive pure local privacy. Our lower bounds allow biased estimators and repeated user participation under the total pure local privacy budget in \eqref{eq:privacy}.

For CDFs, \citet[Theorem~4.1]{LiuHuKong2024} establish uniform and $L^2$ error bounds for the entire CDF curve using randomized threshold comparisons and constrained isotonic estimation under noninteractive pure $\varepsilon$-LDP. Our CDF result determines the exact minimax order of expected integrated squared error over a fixed positive smooth density class, with a matching lower bound allowing finite public blackboard interaction. \citet{HuLiu2026} develop multi-category CDF estimation from censoring with plausible deniability under asymmetric local privacy, also treating symmetric local privacy, and establish uniform consistency and pointwise weak convergence. Their censoring mechanism and inferential results complement our integrated minimax analysis over the specified protocol class.

The rest of the paper is organized as follows. Section~\ref{sec:model} defines the experiment and loss. The sharp rates in Section~\ref{sec:theory} lead to the construction and lower bound in Sections~\ref{sec:estimator}--\ref{sec:lower}, then to CDF, quantile, and transport risks in Section~\ref{sec:cdf}. Remaining proofs are in the appendices.

\section{Statistical model and critical loss}
\label{sec:model}

Let $X_1,\ldots,X_n$ be independent observations from a common unknown density $f$, with user $i$ holding $X_i$. Our goal is to estimate the underlying distribution from reports generated by an $\varepsilon$-locally differentially private protocol.

We allow finite public blackboard protocols in which the next speaker may depend on the public history and users may participate repeatedly. For each $n$, the total number of public messages is bounded by a deterministic finite constant, and a message from user $i$ has a kernel depending only on $X_i$ and the public history. Public randomness is independent of the data and forms part of the transcript $\mathsf Z$, while private randomization is independent across users conditional on the public randomness. All data, message, and randomness spaces are standard Borel, and public messages from the statistician depend only on the history and parameter-independent randomness.

For every user $i$, every pair $x_i,x_i'$ of possible observations, every fixed $x_{-i}$, and every measurable transcript event $B$, require
\begin{equation}
Q\left(B\mid x_i,x_{-i}\right)
\le e^\varepsilon Q\left(B\mid x_i',x_{-i}\right).
\label{eq:privacy}
\end{equation}
Thus $\varepsilon$ is the total pure privacy budget for a user's participation in the complete transcript. The complete-transcript formulation accommodates any budget allocation satisfying \eqref{eq:privacy}; a pathwise sum of conditional message budgets bounded by $\varepsilon$ is one sufficient construction. Write $\mathcal Q_\varepsilon^{\mathrm{BB}}\left(n\right)$ for this class. Noninteractive and fresh-user sequential mechanisms are subclasses. All parameter-independent postprocessing is permitted.

The public local kernel structure, together with \eqref{eq:privacy}, defines the experiment and yields the conditional product structure over users used in Appendix~\ref{app:information}. This factorization supports the information bound for finite public blackboard interaction, including repeated visits to users.

Throughout,
\[
0<\varepsilon\le1,\qquad
N=n\varepsilon^2,\qquad
N_\varepsilon=n\left(e^\varepsilon-1\right)^2.
\]
The quantities $N$ and $N_\varepsilon$ are comparable uniformly in this privacy range. All asymptotic assertions concern $N\to\infty$ with the function-class parameters fixed. Unless specified otherwise, $\left\|\cdot\right\|_2$ is the Euclidean norm for a vector and the Lebesgue $L^2$ norm for a function. Constants in $\lesssim$ and $\asymp$ may depend on the fixed class parameters, but not on $n$ or $\varepsilon$. We use $1/\infty=0$ and natural logarithms, except in explicit truncation rules involving $\log_2$.

We next describe the densities and the resolution structure of the queries.

Let $\T^d=\left[0,1\right]^d$ with opposite faces identified and normalized Lebesgue measure. For $u\in L^2\left(\T^d\right)$, set
\[
\widehat u\left(k\right)=\int_{\T^d}u\left(x\right)e^{-2\pi i k\cdot x}\,dx,
\qquad
\left\|u\right\|_{H^\beta}^2
=\sum_{k\in\mathbb Z^d}\left(1+4\pi^2|k|_2^2\right)^\beta
|\widehat u\left(k\right)|^2.
\]
Fix an integer $d\ge1$, a real $\beta>0$, and $R>1$, and define
\[
\mathcal F_\beta\left(R\right)
=\left\{f\in H^\beta\left(\T^d\right): f\ge0\ \text{a.e.},\
\int_{\T^d}f=1,\ \left\|f\right\|_{H^\beta}\le R\right\}.
\]
The uniform density has norm one. The condition $R>1$ supplies the interior slack needed for the lower bound; when $R=1$, the class consists only of the uniform density. For a density $f$, write $P_f\left(dx\right)=f\left(x\right)\,dx$.

For $\ell=0,1,\ldots$, let
\[
A_\ell=\left\{k\in\mathbb Z^d:2^\ell\le |k|_\infty<2^{\ell+1}\right\},
\qquad D_\ell=|A_\ell|,\qquad \rho_\ell=2^{-d\ell},
\]
and let $P_\ell$ be the orthogonal Fourier projection onto $A_\ell$. Choose one representative from each pair $\left\{k,-k\right\}$ in $A_\ell$, and use both $\sqrt2\cos\left(2\pi k\cdot x\right)$ and $\sqrt2\sin\left(2\pi k\cdot x\right)$. Denote these $D_\ell$ real basis functions by $\psi_{\ell j}$, $1\le j\le D_\ell$. Then
\[
D_\ell=\left(2^{\ell+2}-1\right)^d-\left(2^{\ell+1}-1\right)^d\asymp 2^{d\ell},
\qquad \left\|\psi_{\ell j}\right\|_\infty\le B_0:=\sqrt2.
\]
Write $v_\ell\left(x\right)=\left(\psi_{\ell j}\left(x\right)\right)_{j=1}^{D_\ell}$ and
$\theta_\ell\left(f\right)=\left(\int f\psi_{\ell j}\right)_{j=1}^{D_\ell}$.

The same Fourier blocks define the query class. Its fine index $q$ controls how much variation a test function may allocate across resolutions.

For $1\le q\le\infty$, we equip the periodic Besov space $B^d_{2,q}(\T^d)$ with the following equivalent norm
\begin{equation}
\left\|h\right\|_{\mathcal G_q}
=|\widehat h\left(0\right)|+
\left\|\left(2^{d\ell}\left\|P_\ell h\right\|_2\right)_{\ell\ge0}
\right\|_{\ell_q}.
\label{eq:witness-norm}
\end{equation}
For the standard definition and equivalence across smooth Littlewood-Paley resolutions, see \citet[Definition~3.6 and Theorem~3.7]{BarrazaEtAl2015}; Appendix~\ref{app:fourier} verifies the equivalence of \eqref{eq:witness-norm} with that definition, including both fine-index endpoints. The associated IPM is
\[
d_q\left(\mu,\nu\right)
=\sup_{\left\|h\right\|_{\mathcal G_q}\le1}
\left|\int_{\T^d}h\,d\left(\mu-\nu\right)\right|.
\]
The test functions have continuous bounded representatives, so this definition includes arbitrary probability measures. We also write $d_q\left(u,v\right)$ for mass-one $L^2$ functions, allowing signed functions when constructing an intermediate estimator.

Let $r=q'$ be the conjugate index. For a mass-zero $u\in L^2\left(\T^d\right)$, blockwise Hilbert space duality gives the exact identity
\begin{equation}
d_q\left(u,0\right)
=\left\|\left(\rho_\ell\left\|P_\ell u\right\|_2\right)_{\ell\ge0}
\right\|_{\ell_r}.
\label{eq:block-duality}
\end{equation}
The dual loss therefore aggregates Euclidean errors within each block and uses the $\ell_r$ norm across blocks. For a probability-measure estimate $\widehat P$, define
$\widehat\theta_{\ell j}=\int\psi_{\ell j}\,d\widehat P$. Finite Fourier witnesses give, for every $K$,
\begin{equation}
d_q\left(\widehat P,P_f\right)
\ge
\left\|\left(
\rho_\ell\left\|\widehat\theta_\ell-\theta_\ell\left(f\right)\right\|_2
\right)_{0\le\ell<K}\right\|_{\ell_r}.
\label{eq:finite-duality}
\end{equation}
This inequality supplies the projected loss used in the lower bound for arbitrary probability-measure estimates.

For $p\in\left\{1,2\right\}$, define
\[
\mathcal R_{n,\varepsilon}^{\left(p\right)}\left(q\right)
=\inf_{Q\in\mathcal Q_\varepsilon^{\mathrm{BB}}\left(n\right),\,\widehat P}
\sup_{f\in\mathcal F_\beta\left(R\right)}
\E_f\left[d_q\left(\widehat P,P_f\right)^p\right],
\]
where $\widehat P$ is any measurable probability-measure-valued estimate based on the transcript. Randomized estimates can be included by adding independent postprocessing randomness. The lower proof also applies to arbitrary real estimates of the selected Fourier coefficients under the projected loss in \eqref{eq:finite-duality}.

\section{Sharp minimax rates at critical smoothness}
\label{sec:theory}

We now determine the minimax rates at critical query smoothness across the full Besov fine-index range, with a single noninteractive procedure attaining both first- and second-moment guarantees.

\begin{theorem}[Critical Besov IPM risk]
\label{thm:critical}

Fix $d\ge1$, $\beta>0$, $R>1$, and $1\le q\le\infty$. Uniformly for $0<\varepsilon\le1$ and all sufficiently large $N=n\varepsilon^2$,
\begin{equation}
\mathcal R_{n,\varepsilon}^{(p)}(q)
\asymp N^{-p/2}(\log N)^{p(3/2-1/q)},
\qquad p\in\{1,2\}.
\label{eq:main-rate}
\end{equation}
The comparison constants and the threshold for $N$ depend only on $d,\beta,R,q$.

The upper bounds are attained by the noninteractive mechanism and measurable density estimate in Section~\ref{sec:estimator}. The mechanism and this proper output rule are independent of $q$, and the output belongs to $\mathcal F_\beta\left(R\right)$. The lower bounds hold for all finite public blackboard protocols in $\mathcal Q_\varepsilon^{\mathrm{BB}}\left(n\right)$, with arbitrary bias and repeated access to users subject to the total budget \eqref{eq:privacy}.
\end{theorem}

Equation~\eqref{eq:main-rate} includes both $\E d_q$ and $\E d_q^2$. Appendices~\ref{app:upper} and~\ref{app:lower} prove the bounds; equivalent Besov norms change only the constants.

\begin{table}[!ht]
\centering
\begin{tabular}{@{}ccll@{}}
\toprule
Test index $q$ & Error aggregation & Expected IPM risk & Squared IPM risk \\
\midrule
$1$ & $\ell_\infty$ over blocks & $\sqrt{\log N/N}$ & $\left(\log N\right)/N$ \\
$2$ & $\ell_2$ over blocks & $\left(\log N\right)/\sqrt N$ & $\left(\log N\right)^2/N$ \\
$\infty$ & $\ell_1$ over blocks & $\left(\log N\right)^{3/2}/\sqrt N$ & $\left(\log N\right)^3/N$ \\
\bottomrule
\end{tabular}
\caption{Endpoint and Sobolev cases of Theorem~\ref{thm:critical}. Aggregation refers to the weighted block errors in \eqref{eq:block-duality}. Corollary~\ref{cor:smoothness-free} attains the same risk orders without using smoothness, with output in $\mathcal F_0\left(R\right)$.}
\label{tab:endpoints}
\end{table}

At $q=2$, the test norm is equivalent to periodic $H^d$. To compare with~\citet{AlbertEtAl2026}, write $\mathcal A_\beta\left(R_A\right)$ for their density ball of radius $R_A$ and $d_{A,R_D}$ for their IPM with discriminator radius $R_D$, as defined in Lemma~\ref{lem:albert-classes}. Their squared norm of the uniform density is $d$, so the nondegenerate comparison requires fixed $\beta>0$, $R_A>\sqrt d$, and $R_D>0$. The lemma gives radii $R_-,R_+>1$ in \eqref{eq:albert-radii} such that
\[
 \mathcal F_\beta\left(R_-\right)
 \subseteq\mathcal A_\beta\left(R_A\right)
 \subseteq\mathcal F_\beta\left(R_+\right),
\]
and \eqref{eq:albert-loss-comparison} compares $d_{A,R_D}$ with $R_Dd_2$ by fixed constants. Applying the lower bound on the inner ball and the upper bound on the outer ball transfers the rate to their class. The upper estimator takes values in the outer Sobolev ball; no membership in $\mathcal A_\beta\left(R_A\right)$ is asserted. At $R_A=\sqrt d$ the density class contains only the uniform density, and below that radius it is empty. These comparisons close the gap between $N^{-1/2}$ and $N^{-1/2}\left(\log N\right)^{2d}$. The cited lower bound assumes integer smoothness.

Without privacy, the empirical measure already satisfies $\sup_f\E_f d_q\left(P_n,P_f\right)^2\le C_d/n$ for every $q$: each Fourier block has mean squared coefficient error at most $B_0^2D_\ell/n$, and $\sum_\ell\rho_\ell\sqrt{D_\ell}<\infty$. Appendix~\ref{app:fourier} gives the calculation. Thus privacy distinguishes these losses through an additional logarithmic cost.

The critical identity $D_\ell\rho_\ell\asymp1$ makes each retained resolution contribute the same order to the ellipsoid normalization. With $K$ levels, the maximum weighted block error has first and second moments of order $\sqrt{K/N}$ and $K/N$. The aggregation inequalities
\begin{equation}
K^{1/r-1}\sum_{\ell<K}t_\ell
\le\left\|\left(t_\ell\right)_{\ell<K}\right\|_{\ell_r}
\le K^{1/r}\max_{\ell<K}t_\ell,
\qquad t_\ell\ge0,
\label{eq:aggregation}
\end{equation}
then explain all intermediate indices. The upper and lower arguments control the right and left sides of \eqref{eq:aggregation}, respectively. With $K\asymp\log N$, this gives $1/2+1/r=3/2-1/q$.

The bias $2^{-\left(\beta+d\right)K}$ determines the constant relating $K$ to $\log N$. Fixed smoothness changes this constant, leaving the critical logarithmic exponent unchanged.

\begin{corollary}[Known radius and unknown smoothness]
\label{cor:smoothness-free}
\label{rem:smoothness-free}
Let $d\ge1$, $R>1$ be known, and
$\mathcal F_0(R)=\{f\ge0:\int_{\T^d}f=1,\ \|f\|_2\le R\}$.
For $0<\varepsilon\le1$ and $N\ge2$, use
\eqref{eq:ellipsoid}--\eqref{eq:decoder} with
$K_0=\lceil\log_2N/(2d)\rceil$, and \eqref{eq:selector} with
tolerance $N^{-1/2}$ and a fixed countable $L^2$-dense family in
$\mathcal F_0(R)$. The measurable output $\widehat f^{(0)}$ satisfies,
for $C=C(d,R)$ and all $1\le q\le\infty$, $p\in\{1,2\}$,
\begin{equation}
 \sup_{f\in\mathcal F_0(R)}\E_f d_q(\widehat f^{(0)},f)^p
 \le C N^{-p/2}(\log N)^{p(3/2-1/q)}.
 \label{eq:smoothness-free-risk}
\end{equation}
This procedure uses $n,\varepsilon,d,R$ but neither $\beta$ nor $q$,
and is minimax rate-optimal on every fixed $\mathcal F_\beta(R)$,
$\beta\ge0$, with output in $\mathcal F_0(R)$.
The lower-bound constants and thresholds may depend on $\beta$.
Each user sends one bit, with release dimension $m_{K_0}\asymp\sqrt N$
instead of $m_K\asymp N^{d/(2(\beta+d))}$ in \eqref{eq:cutoff}.
\end{corollary}

The proof is in Appendix~\ref{app:smoothness-free-proof}.
Proposition~\ref{prop:f0-finite} gives finite-dimensional convex
postprocessing of the same release, with output in $\mathcal F_0(R)$
and an explicit optimization tolerance preserving these rates.
\section{An explicit noninteractive release}
\label{sec:estimator}

The construction balances errors across resolutions while preserving their different concentration strengths. Ellipsoid scaling gives each block the same weighted error scale; finer blocks have more coordinates but smaller directional noise. Retaining their dimensions in the tail bounds controls the maximum over all blocks. We implement this scaling with a random direction and one private bit per user, then select a density that supports every $q$.

For $N>1$, retain the levels $0,\ldots,K-1$, where
\begin{equation}
K=\left\lceil\frac{\log_2 N}{2\left(\beta+d\right)}\right\rceil.
\label{eq:cutoff}
\end{equation}
Let $m_K=\sum_{\ell<K}D_\ell$ and $S_K=\sum_{\ell<K}D_\ell\rho_\ell$. For $x\in\T^d$, form the weighted feature vector and block-diagonal matrix
\begin{equation}
z\left(x\right)=\left(\rho_\ell v_\ell\left(x\right)\right)_{\ell<K},
\qquad
A\big|_\ell=B_0\sqrt{S_K\rho_\ell}\,I_{D_\ell},
\qquad t\left(x\right)=A^{-1}z\left(x\right).
\label{eq:ellipsoid}
\end{equation}
Then $\left\|t\left(x\right)\right\|_2\le1$ for all $x$, and $S_K\asymp K$.

Put $\kappa=\tanh\left(\varepsilon/2\right)$ and let
\[
c_m=\E|U_1|
=\frac{\Gamma\left(m/2\right)}{\sqrt\pi\,\Gamma\left(\left(m+1\right)/2\right)},
\qquad U\sim\operatorname{Unif}\left(\mathbb S^{m-1}\right).
\]
The release is as follows.
\begin{enumerate}
\item Independently for each user, generate a public direction $U_i$ uniformly on $\mathbb S^{m_K-1}$, independently of all observations.
\item User $i$ sends one bit $B_i\in\left\{-1,1\right\}$, with
\begin{equation}
\Prob(B_i=b\mid U_i,X_i=x)
=\frac12\left[1+b\kappa\|t(x)\|_2
\sign\langle U_i,t(x)\rangle\right],
\qquad b\in\{-1,1\}.
\label{eq:bit}
\end{equation}
Use $\sign\left(0\right)=0$, which also specifies a fair bit when $t\left(x\right)=0$.
\item The statistician decodes
\begin{equation}
W_i=\frac{AB_iU_i}{\kappa c_{m_K}},
\quad \overline W=\frac1n\sum_{i=1}^nW_i,
\quad Y(x)=1+\sum_{\ell<K}\rho_\ell^{-1}
\langle\overline W_\ell,v_\ell(x)\rangle.
\label{eq:decoder}
\end{equation}
\end{enumerate}
The decoder $W_i$ estimates the weighted features, and $\rho_\ell^{-1}\overline W_\ell$ estimates the density coefficients. The mechanism uses independent real-valued public directions, with each user communicating one private response bit.

Conditionally on any direction, the bit probabilities in \eqref{eq:bit} lie between $\left(1-\kappa\right)/2$ and $\left(1+\kappa\right)/2$, whose ratio is $e^\varepsilon$. Hence the mechanism is pure $\varepsilon$-LDP. For $U$ uniform on $\mathbb S^{m_K-1}$, spherical symmetry gives
\[
\E_U\left[U\sign\langle U,t\rangle\right]
=c_{m_K}\frac{t}{\|t\|_2},\qquad t\ne0.
\]
Consequently, the radial factor $\|t(x)\|_2$ in \eqref{eq:bit} gives $\E[W_i\mid X_i=x]=z(x)$, including at $t(x)=0$, so $Y$ has mass one and is unbiased for the truncated Fourier series of $f$. The construction is a public-direction implementation of established spherical local randomization, after the ellipsoid transformation~\citep{DuchiJordanWainwright2018,NikolovTang2024}.

To turn the mass-one Fourier estimate $Y$ into a density using the same reports, fix once and for all a countable $L^2$-dense family $\left(g_m\right)_{m\ge1}$ in $\mathcal F_\beta\left(R\right)$, and let $\tau=N^{-1/2}$. Define the approximate minimum-distance selector
\begin{equation}
I\left(Y\right)=\inf_{m\ge1}d_1\left(Y,g_m\right),\qquad
\widehat f=g_{m\left(Y\right)},\qquad
m\left(Y\right)=\min\left\{m:d_1\left(Y,g_m\right)\le I\left(Y\right)+\tau\right\}.
\label{eq:selector}
\end{equation}
The countable selector defines a measurable statistical rule. Proposition~\ref{prop:f0-finite} gives a finite-dimensional alternative with output in the known $L^2$ ball.
This is a measurable rule with $\widehat f\in\mathcal F_\beta\left(R\right)$ and
\[
d_1\left(\widehat f,f\right)\le2d_1\left(Y,f\right)+\tau.
\]
The block-maximum distance $d_1$ defines a common output rule for every fine index $q$. Since both $\widehat f$ and $f$ lie in the Sobolev ball, their high-frequency tails are controlled, and applying \eqref{eq:aggregation} to the first $K$ blocks gives the upper bound for every $q$. Appendix~\ref{app:proper-output} establishes the risk bounds and measurability of this approximate minimum-distance rule.

The remaining step controls the largest weighted block error jointly across resolutions. Write
\[
E_\ell=\overline W_\ell-\rho_\ell\theta_\ell\left(f\right),
\qquad W_K^*=\max_{\ell<K}\left\|E_\ell\right\|_2.
\]
The deterministic identity
$|\left\langle a,W_i\right\rangle|=|\left\langle A^\top a,U_i\right\rangle|/\left(\kappa c_{m_K}\right)$
allows the spherical moment bound to be used even though the bit depends on the direction. A net argument in a block of dimension $D_\ell$ gives, for a sufficiently large fixed $A_0$ and $u\ge1$,
\begin{equation}
\Prob_f\left(\left\|E_\ell\right\|_2>A_0\sqrt{K/N}\,u\right)
\le2e^{-cD_\ell u^2},\qquad
\Prob_f\left(W_K^*>A_0\sqrt{K/N}\,u\right)
\le Ce^{-c'u^2}.
\label{eq:block-tail}
\end{equation}
The second bound follows by summing the first over geometrically increasing $D_\ell$, uniformly in $K$. In particular, $\E_f\left(W_K^*\right)^p\lesssim\left(K/N\right)^{p/2}$ for $p=1,2$. Retaining the dimensions in \eqref{eq:block-tail} is what avoids an additional iterated logarithm in the maximum-block loss.

The cutoff makes the Sobolev tail at most a constant depending only on $d,\beta,R$ times $N^{-1/2}$. Combining the selector inequality, the aggregation bound \eqref{eq:aggregation}, and the first and second moments from \eqref{eq:block-tail} gives the upper rates in Theorem~\ref{thm:critical}. Appendix~\ref{app:proper-output} contains the complete calculation.

For the known-radius construction in Corollary~\ref{cor:smoothness-free}, the tail is controlled without using $\beta$. Set $K_0=\lceil\log_2N/(2d)\rceil$ and $u=\widehat f^{(0)}-f$, where $\widehat f^{(0)},f\in\mathcal F_0(R)$, so $\|u\|_2\le2R$. Since $\|P_\ell u\|_2\le\|u\|_2$, summing the geometric weights $\rho_\ell=2^{-d\ell}$ bounds the omitted blocks in every $\ell_r$ norm, $1\le r\le\infty$, by $\|u\|_2 2^{-dK_0}/(1-2^{-d})\le4RN^{-1/2}$; see \eqref{eq:smoothness-free-tail}. The retained blocks obey the same first and second moment bounds with $K_0$ in place of $K$. Combining this tail bound with the selector inequality and \eqref{eq:aggregation}, with $r=q'$, yields the rates for every $q$ using only the known $L^2$ radius. Appendix~\ref{app:smoothness-free-proof} gives the complete argument, and the finite-dimensional implementation is given in Proposition~\ref{prop:f0-finite}.

\section{Interaction cannot improve the minimax rates}
\label{sec:lower}

The noninteractive rates in Theorem~\ref{thm:critical} remain minimax optimal over finite public blackboard protocols, including repeated visits under each user's total privacy budget. We establish this by distributing a prior across $K$ resolutions, with coordinate widths chosen so that each level contributes the same prior information while balancing the loss weights. As $N$ grows, the entire family approaches the uniform density and remains inside the fixed smooth class, preserving the joint difficulty of the retained resolutions. Appendices~\ref{app:information} and~\ref{app:lower} provide the information inequality and proof details.

Use $K=\left\lfloor\gamma\log_2N_\varepsilon\right\rfloor$ levels, where
$0<\gamma<\min\left\{\left(3d\right)^{-1},\left(2\beta+2d\right)^{-1}\right\}$, and consider
\[
f_\theta=1+\sum_{\ell<K}\sum_{j=1}^{D_\ell}\theta_{\ell j}\psi_{\ell j},
\qquad |\theta_{\ell j}|\le a_\ell,
\qquad a_\ell^2=\frac{c_0KD_\ell}{N_\varepsilon},
\]
with a fixed $c_0>0$. These widths make both $\left\|f_\theta-1\right\|_\infty$ and
$\left\|f_\theta\right\|_{H^\beta}^2-1$ tend to zero uniformly over the rectangle.
Consequently its entire support eventually lies in $\mathcal F_\beta\left(R\right)$,
with $1/2\le f_\theta\le3/2$. Let $\Pi$ be the product of cosine-squared
priors on the coordinates. Each coordinate has prior Fisher information $\pi^2/a_\ell^2$,
so the total prior information equals $\pi^2N_\varepsilon/c_0$.

This prior allocates information across the coefficients. The transcript has a matching global budget: the unprivatized score covariance is bounded by $2I$. The private information
contraction in Lemma~\ref{lem:fisher} gives
$\operatorname{tr}I_{\mathsf Z}\left(\theta\right)\le2N_\varepsilon$ for every permitted
protocol. In the public blackboard model, the observations remain
conditionally independent given the transcript: each message updates
only its speaker's posterior factor. This remains true when users are
revisited. It is this factorization, together with each user's total
privacy constraint, that permits the same information budget under
interaction.

To convert this information budget into the required loss, clip each estimated Fourier coefficient to its prior interval and write
$e_\ell$ for the resulting error vector in block $\ell$. Clipping reduces
each coefficient error and gives $|e_{\ell j}|\le2a_\ell$. The coordinate van Trees
inequality in Lemma~\ref{lem:van-trees} then yields the first-moment bound
\eqref{eq:coordinate-first-moment}, using
$\E|e_{\ell j}|\ge\E e_{\ell j}^2/\left(2a_\ell\right)$.

Put $u_\ell=\rho_\ell/\left(2a_\ell\sqrt{D_\ell}\right)$. Weighted Cauchy-Schwarz and the two information budgets give
\[
\E_\Pi\sum_{\ell<K}\rho_\ell\left\|e_\ell\right\|_2
\ge \frac{\left(\sum_{\ell<K}D_\ell\sqrt{u_\ell}\right)^2}{\left(2+\pi^2/c_0\right)N_\varepsilon}
\gtrsim \frac{K^{3/2}}{\sqrt{N_\varepsilon}},
\]
because
\[
D_\ell\sqrt{u_\ell}
=2^{-1/2}\left(\frac{N_\varepsilon}{c_0K}\right)^{1/4}
\sqrt{D_\ell\rho_\ell}\asymp\left(N_\varepsilon/K\right)^{1/4}.
\]
Thus all $K$ resolutions contribute equally; see \eqref{eq:lower-weighted-cauchy}.
 Applying the left inequality in
\eqref{eq:aggregation} gives expected loss at least a constant times
$K^{1/2+1/r}/\sqrt{N_\varepsilon}$ for every $r=q'$, including both
endpoints. Jensen's inequality gives the corresponding squared-risk
bound. Since $K\asymp\log N$ and $N_\varepsilon\asymp N$, these are the
lower rates in Theorem~\ref{thm:critical}.

\section{CDFs, quantiles, and transport}
\label{sec:cdf}

Estimating an entire CDF from one set of locally private reports requires accurate threshold probabilities across $\left[0,1\right]$. While binary randomized response estimates one prespecified threshold with squared risk $O\left(N^{-1}\right)$, the minimax integrated squared CDF error over fixed positive smooth density classes is $\left(\log N\right)^2/N$. We establish this rate through a shared information-budget lower bound and hierarchical range estimation, then use monotone projection and inversion to obtain sharp quantile, Wasserstein, and transport map risks.

We use H\"older regularity on the interval $\left[0,1\right]$ to define the density class. For noninteger $s=m+\alpha>0$, $m\ge0$, $0<\alpha<1$, use
\[
\left\|u\right\|_{C^s}
=\sum_{j=0}^m\left\|u^{\left(j\right)}\right\|_\infty
+\sup_{x\ne y}\frac{|u^{\left(m\right)}\left(x\right)-u^{\left(m\right)}\left(y\right)|}{|x-y|^\alpha}.
\]
For integer $s=m\ge1$, use $\sum_{j=0}^m\left\|u^{\left(j\right)}\right\|_\infty$ on functions with continuous derivatives through order $m$. In either convention $\left\|1\right\|_{C^s}=1$. Fix $M>1$ and $0<a<1<b$, and set
\[
\mathcal H_s\left(M,a,b\right)
=\left\{f\in C^s\left(\left[0,1\right]\right):\ \textstyle\int_0^1 f=1,
\ \left\|f\right\|_{C^s}\le M,\ a\le f\le b\right\}.
\]
Write $F_f\left(x\right)=\int_0^x f\left(t\right)\,dt$. The public local model of Section~\ref{sec:model} is used with one interval-valued observation per user.

\begin{theorem}[Integrated CDF error]
\label{thm:cdf}
For fixed $s>0$, $M>1$, and $0<a<1<b$, and uniformly for $0<\varepsilon\le1$ and all sufficiently large $N=n\varepsilon^2$,
\begin{equation}
\inf_{Q\in\mathcal Q_\varepsilon^{\mathrm{BB}}\left(n\right),\,\widehat F}
\sup_{f\in\mathcal H_s\left(M,a,b\right)}
\E_f\int_0^1|\widehat F\left(x\right)-F_f\left(x\right)|^2\,dx
\asymp\frac{\left(\log N\right)^2}{N}.
\label{eq:cdf-risk}
\end{equation}
A noninteractive procedure attains the upper bound with output equal to the CDF of a density bounded between $a$ and $b$, while the lower bound holds for arbitrary measurable function estimates under finite public blackboard interaction, including repeated visits to users.
\end{theorem}

The squared logarithm persists at every fixed finite smoothness: an anchored cosine family gives the lower bound, and a hierarchical range estimator tuned using $N,a,b$ attains the upper bound. Smoothness determines the lower family and the asymptotic threshold.

The interval CDF is linked to the periodic negative Sobolev loss through an anchoring constant. Specifically, if $\mathcal P\left(f-1\right)$ is the zero-mean periodic primitive of $f-1$, then
\begin{equation}
F_f\left(x\right)=x+\mathcal P\left(f-1\right)\left(x\right)+\frac12-\E_f X.
\label{eq:interval-anchor}
\end{equation}
The constant term in \eqref{eq:interval-anchor} fixes the interval boundary values and contributes separately from the zero-mean periodic primitive. For two interval densities, \eqref{eq:cdf-periodic-decomposition} separates the periodic primitive loss from the squared difference of their means. Theorem~\ref{thm:cdf} uses the interval class $\mathcal H_s(M,a,b)$ and the direct constructions below: the periodic $d=1,q=2$ result alone does not include this anchoring term or identify the two density classes. We estimate interval CDFs directly through range queries and construct the lower bound from cosine perturbations of the uniform density, whose anchored CDFs are sine series plus $x$. Pairing a function estimate $H-x$ with the corresponding sine recovers an estimate of each cosine coefficient, as detailed in \eqref{eq:cdf-coefficient-lower} in Appendix~\ref{app:cdf-lower}.

For the cosine family, the integrated CDF loss weights coefficient $k$
by $w_k=\left(2\pi k\right)^{-2}$. Put $S_J=\sum_{k\le J}\sqrt{w_k}\asymp\log J$
and choose squared prior widths $a_k^2=tS_J/\left(N_\varepsilon\sqrt{w_k}\right)$
with fixed $t>0$. The total prior information is then of order
$N_\varepsilon$. Coordinate van Trees and weighted Cauchy-Schwarz give
Bayes risk at least $cS_J^2/N_\varepsilon$. Meanwhile the whole rectangle
satisfies
$\left\|f_\theta-1\right\|_{C^s}\lesssim\sqrt{\log J/N_\varepsilon}\,J^{s+3/2}$.
Taking $J=\left\lfloor N_\varepsilon^\eta\right\rfloor$ with
$0<\eta<\left(2s+3\right)^{-1}$ keeps the family inside the fixed positive smooth
class and gives the lower rate $\left(\log N\right)^2/N$. The information budget
again covers the full public transcript. Appendix~\ref{app:cdf-lower}
justifies the coefficient reduction for arbitrary function estimates.

To attain the upper bound, we use hierarchical range estimation \citep{CormodeKulkarniSrivastava2019}, assigning users in advance to disjoint groups of nearly equal size, one for each retained dyadic level. Each user spends the full privacy budget on an optimized unary encoding of the cell containing the observation, and a grid-prefix probability is reconstructed from at most one cell per level. Independence between the user groups makes the prefix variance the sum of its levelwise variances, even when cell estimates within a level are dependent. This allocation yields the squared logarithm: with $h=\left\lceil\tfrac12\log_2N\right\rceil$ and $D=2^h$, the decoded step function $\widetilde F$ satisfies
\[
\E_f\left\|\widetilde F-F_f\right\|_2^2
\le \frac{Ch^2}{N}+\frac{b^2}{D^2}.
\]
Let $\mathcal C_{a,b}$ consist of the continuous functions with
$G\left(0\right)=0$, $G\left(1\right)=1$, and
\begin{equation}
 a\left(y-x\right)\le G\left(y\right)-G\left(x\right)\le b\left(y-x\right),\qquad 0\le x<y\le1.
\label{eq:cdf-increments}
\end{equation}
Viewed as a subset of $L^2\left[0,1\right]$, this set is nonempty, closed,
and convex. Define its Hilbert projection by
\begin{equation}
 \widehat F=\operatorname*{argmin}_{G\in\mathcal C_{a,b}}
 \left\|G-\widetilde F\right\|_2^2.
\label{eq:cdf-projection}
\end{equation}
The projection is unique and measurable. Since $F_f\in\mathcal C_{a,b}$,
\begin{equation}
 \left\|\widehat F-F_f\right\|_2^2\le\left\|\widetilde F-F_f\right\|_2^2.
\label{eq:cdf-projection-risk}
\end{equation}
Thus $\widehat F$ attains \eqref{eq:cdf-risk} and is the CDF of a density bounded between $a$ and $b$, with the report construction proved in Appendix~\ref{app:cdf-hierarchy} and the projection properties in Appendix~\ref{app:cdf-projection}. Proposition~\ref{prop:cdf-projection-gap} gives a finite-dimensional version on the same $D$-cell grid, using linear endpoint and slope constraints and a quadratic objective. For a feasible output with objective error $\tau$, the expected squared CDF error is at most three times the sum of $\tau$, the release risk, and $4b^2/D^2$, so an optimization error of order $\left(\log N\right)^2/N$ preserves the rate. The slope bounds also transfer this control to quantiles and transport maps.

\subsection{Quantiles and distribution estimation}
The same hierarchical reports estimate the full quantile function and the underlying distribution through projection onto the set $\mathcal D$ of all $\left[0,1\right]$-valued nondecreasing $L^2$ functions. The projected CDF defines a probability law, allowing atoms, and its generalized inverse estimates $q_f=F_f^{-1}$. This postprocessing uses only the reports, $n$, and $\varepsilon$.

\begin{corollary}
\label{cor:quantile}
Under the fixed class parameters and privacy regime of
Theorem~\ref{thm:cdf},
\[
 \inf_{Q\in\mathcal Q_\varepsilon^{\mathrm{BB}}\left(n\right),\,\widehat q}
 \sup_{f\in\mathcal H_s\left(M,a,b\right)}
 \E_f\left\|\widehat q-q_f\right\|_2^2
 \asymp \frac{\left(\log N\right)^2}{N},
\]
and
\[
 \inf_{Q\in\mathcal Q_\varepsilon^{\mathrm{BB}}\left(n\right),\,\widehat P}
 \sup_{f\in\mathcal H_s\left(M,a,b\right)}
 \E_f W_2^2\left(\widehat P,P_f\right)
 \asymp \frac{\left(\log N\right)^2}{N}.
\]
The first infimum allows arbitrary measurable function estimates;
the second allows arbitrary probability measures on $\left[0,1\right]$.
A single noninteractive release, tuned using only $n,\varepsilon$, attains both upper bounds.
\end{corollary}

For any probability measure $P$ on $\left[0,1\right]$ with CDF $G$, the target bounds
$a\le f\le b$ give
\begin{equation}
 aW_2^2\left(P_f,P\right)\le\left\|F_f-G\right\|_2^2\le bW_2^2\left(P_f,P\right).
\label{eq:cdf-wasserstein-comparison}
\end{equation}
Appendix~\ref{app:cdf-wasserstein-proof} proves this comparison, including atomic $P$ and endpoint atoms.

The set $\mathcal D$ is nonempty, convex, and closed: Helly's
selection theorem and dominated convergence identify the limit of every
$L^2$-convergent sequence with a member of $\mathcal D$.
Project the hierarchical step release onto this set,
$\widehat G=\Pi_{\mathcal D}\widetilde F$.
For the hierarchical step release, this projection is the finite
bounded isotonic least-squares problem
\begin{equation}
 \widehat{\boldsymbol z}
 =\operatorname*{argmin}_{0\le z_0\le\cdots\le z_{D-1}\le1}
       \frac1D\sum_{j=0}^{D-1}(z_j-\widehat F_j)^2,
 \qquad
 \widehat G=\sum_{j=0}^{D-1}\widehat z_j
       \mathbf1_{[j/D,(j+1)/D)}.
 \label{eq:cdf-finite-isotonic}
\end{equation}
Cell averaging reduces the $L^2$ projection to this problem;
Appendix~\ref{app:finite-isotonic} proves the reduction and continuity
of its unique solution.
Each member determines a probability law by its right-continuous
representative on $\left[0,1\right)$ and CDF value one at $1$. Denote the resulting
law by $\widehat P$ and output its generalized quantile $\widehat q$.
Projection contracts CDF error because $F_f\in\mathcal D$, so
\eqref{eq:cdf-wasserstein-comparison} gives
\[
 \E_f\left\|\widehat q-q_f\right\|_2^2
 =\E_fW_2^2\left(\widehat P,P_f\right)
 \le a^{-1}\E_f\left\|\widehat G-F_f\right\|_2^2
 \le a^{-1}\E_f\left\|\widetilde F-F_f\right\|_2^2.
\]
Equation~\eqref{eq:cdf-hierarchy-bound} proves the rate. The density bound
$a$ enters this risk comparison, while the procedure itself uses only
$n,\varepsilon$. Its measurability also follows directly: the projection
is nonexpansive in $L^2$, and for laws corresponding to $G,H\in\mathcal D$,
\[
 W_2^2\left(P_G,P_H\right)\le W_1\left(P_G,P_H\right)
 =\left\|G-H\right\|_1\le\left\|G-H\right\|_2.
\]
The first inequality uses $|x-y|^2\le|x-y|$ on $\left[0,1\right]$ and the middle
identity follows by integrating the region between the two graphs.
Thus the induced law and its $L^2$ quantile depend continuously on the
projected CDF class.

Appendix~\ref{app:quantile} proves the lower bounds, while the slope-constrained output in \eqref{eq:cdf-projection} supplies the invertible CDFs used for transport maps below.

\subsection{Transport between two smooth distributions}

Let $X_1,\ldots,X_n\sim f$ and $Y_1,\ldots,Y_n\sim g$ be independent labeled samples held by $2n$ distinct users, each holding one observation and having total privacy budget $\varepsilon$. A public protocol may involve all $2n$ users, including repeated visits. The increasing map $T_{f,g}=F_g^{-1}\circ F_f$ is the one-dimensional optimal map for quadratic transport cost, and we measure estimation error by
$\left\|\widehat T-T_{f,g}\right\|_{L^2\left(P_f\right)}^2=\int_0^1|\widehat T\left(x\right)-T_{f,g}\left(x\right)|^2 f\left(x\right)\,dx$.

Fix $\Gamma>1$ and $0<L_1<1<L_2$. Let $\mathcal M_s\left(M,\Gamma,L_1,L_2\right)$ consist of pairs $\left(f,g\right)$ of $C^s$ densities with norms at most $M$, satisfying $\Gamma^{-1}<f,g<\Gamma$, for which there is a $C^2$ convex potential $\phi$ with
\[
\phi'=T_{f,g},\qquad L_1\le\phi''\le L_2,\qquad
\phi'\left(\left[0,1\right]\right)\subseteq\left[0,1\right].
\]
This is the one-dimensional class in Assumption~2 of \citet{LalanneEtAl2026}, with fixed interior slack at the uniform pair. Their local privacy lower bound in Table~\ref{tab:closest-rates} assumes sufficiently small privacy budgets.

\begin{corollary}
\label{cor:transport}
Fix $s>0$, $M>1$, $\Gamma>1$, and $0<L_1<1<L_2$. With equal sample sizes $n$ and $N=n\varepsilon^2$, uniformly for $0<\varepsilon\le1$ and sufficiently large $N$,
\begin{equation}
\inf_{Q\in\mathcal Q_\varepsilon^{\mathrm{BB}}\left(2n\right),\,\widehat T}
\sup_{\left(f,g\right)\in\mathcal M_s\left(M,\Gamma,L_1,L_2\right)}
\E_{f,g}\left\|\widehat T-T_{f,g}\right\|_{L^2\left(P_f\right)}^2
\asymp\frac{\left(\log N\right)^2}{N}.
\label{eq:transport-risk}
\end{equation}
Here the protocol class is applied to the two independent labeled samples. A noninteractive upper procedure estimates the two CDFs separately by Theorem~\ref{thm:cdf} and sets $\widehat T=G_Y^{-1}\circ G_X$, using $a=\Gamma^{-1}$ and $b=\Gamma$. The lower holds for arbitrary measurable map estimates.
\end{corollary}

The corollary converts integrated distributional accuracy into map accuracy. For any $G_X,G_Y\in\mathcal C_{a,b}$ and truths $a\le f,g\le b$,
\begin{equation}
\left\|G_Y^{-1}\circ G_X-T_{f,g}\right\|_{L^2\left(P_f\right)}^2
\le\frac{2b}{a^2}
\left(\left\|G_X-F_f\right\|_2^2+\left\|G_Y-F_g\right\|_2^2\right).
\label{eq:transport-stability}
\end{equation}
Appendix~\ref{app:transport-proof} proves the stability bound.

For the lower bound in \eqref{eq:transport-risk}, fix $g=1$. Then $T_{f,1}=F_f$ and $\phi''=f$. The positive cosine prior used for \eqref{eq:cdf-risk} lies in the required potential class because $L_1<1<L_2$. The target observations have zero score for that submodel, so the information bound counts the $n$ informative source users even if all $2n$ users participate repeatedly. Appendix~\ref{app:cdf-lower} gives the full class verification and the weighted-loss comparison that complete the transport lower bound.

\section{Discussion}

Local privacy imposes a shared multiscale information budget whose interaction with the query loss determines the exact logarithmic rates. Dimension-matched priors and dimension-sensitive concentration give matching lower and upper bounds, including first-moment risk and public interaction. A common density output supports every fine index, and a known $L^2$ radius suffices for risk adaptation. The same budget argument yields sharp CDF risk, with quantile, Wasserstein, and transport consequences.

Proposition~\ref{prop:f0-finite} gives finite-dimensional proper output in the known $L^2$ ball, with an explicit tolerance for optimization error. Finite-dimensional output constrained to the original Sobolev ball and confidence sets for smooth expectation queries and CDFs remain natural next directions. Nonnegative trigonometric approximation offers a route to the former, while dimension-sensitive concentration provides a starting point for joint coverage guarantees.
\label{end:main}

\clearpage
\appendix

\section{Fourier blocks, Besov duality, and bias}
\label{app:fourier}

The function-space facts below identify the norm in \eqref{eq:witness-norm} with the usual periodic Besov norm and support both risk bounds.

\subsection{Equivalence with a smooth resolution}

Let $\left(\Delta_j\right)_{j\ge0}$ be a smooth inhomogeneous periodic Littlewood-Paley resolution. Thus its Fourier multipliers form a partition of unity, the low block contains a fixed neighborhood of zero, and the other multipliers are supported in dyadic annuli with uniformly bounded size and overlap. The standard norm is
\[
\left\|\left(2^{dj}\left\|\Delta_jh\right\|_2\right)_{j\ge0}\right\|_{\ell_q}.
\]
For the conventional smooth-resolution definition and independence of the resolution, see~\citep[Definition~3.6 and Theorem~3.7]{BarrazaEtAl2015}.

Away from a fixed finite set of low frequencies, $\Delta_j$ can meet $A_\ell$ only if $|j-\ell|\le L_d$ for a constant $L_d$. Plancherel and the bounded multipliers give
\[
2^{dj}\left\|\Delta_jh\right\|_2
\le C\sum_{|\ell-j|\le L_d}2^{d\ell}\left\|P_\ell h\right\|_2.
\]
Conversely, use the partition of unity and the fact that $P_\ell$ is an $L^2$ contraction to obtain
\[
2^{d\ell}\left\|P_\ell h\right\|_2
\le C\sum_{|j-\ell|\le L_d}2^{dj}\left\|\Delta_jh\right\|_2.
\]
The finitely many low terms in these inequalities are augmented by $|\widehat h\left(0\right)|$ and fixed neighboring blocks. Taking an $\ell_q$ norm and using convolution by a fixed finite sequence proves equivalence in both directions for every $1\le q\le\infty$. In particular, no equivalence constant depends on a truncation level. For $q=2$, Plancherel also shows equivalence with $H^d$ directly.

\subsection{Continuous representatives and duality}

For a function with $\left\|h\right\|_{\mathcal G_q}\le1$, every scaled block norm is at most one. The real basis bound gives
\[
\left\|P_\ell h\right\|_\infty
\le B_0\sqrt{D_\ell}\left\|P_\ell h\right\|_2
\le C2^{-d\ell/2}.
\]
The Fourier block series therefore converges uniformly to a continuous representative, uniformly over the test ball. All integrals against finite measures are consequently defined.

For mass-zero $u\in L^2$, the zero-frequency term pairs to zero, and Cauchy-Schwarz in each block followed by sequence H\"older gives
\[
\left|\int hu\right|
\le\sum_{\ell\ge0}\left\|P_\ell h\right\|_2\left\|P_\ell u\right\|_2
\le\left\|\left(\rho_\ell\left\|P_\ell u\right\|_2\right)_\ell\right\|_{\ell_r}
\quad\text{if }\left\|h\right\|_{\mathcal G_q}\le1.
\]
For the reverse inequality, restrict to the first $K$ blocks. Choose $P_\ell h$ parallel to $P_\ell u$ and give it norm $\rho_\ell b_\ell$, where $b_\ell\ge0$ and $\left\|\left(b_\ell\right)_{\ell<K}\right\|_{\ell_q}\le1$. Ordinary finite-dimensional duality yields the $\ell_r$ norm of the first $K$ weighted error blocks. Letting $K\to\infty$ proves \eqref{eq:block-duality}, including $q=1$ and $q=\infty$. The endpoint argument follows by taking limits of the finite-dimensional test functions.

For a finite signed measure of total mass zero, the same finite witnesses give \eqref{eq:finite-duality}. Uniform convergence of the test expansions, followed by rational approximation of the coefficients of a finite witness, also expresses the IPM on probability measures as a countable supremum of measurable finite integrals. Hence the losses used in the risk definitions are measurable.

For later use, \eqref{eq:block-duality} implies the Lipschitz bound
\begin{equation}
d_q\left(u,0\right)
\le\sum_{\ell\ge0}\rho_\ell\left\|P_\ell u\right\|_2
\le\left(\sum_{\ell\ge0}\rho_\ell^2\right)^{1/2}\left\|u\right\|_2
\le C_d\left\|u\right\|_2.
\label{eq:l2-continuity}
\end{equation}

For $P_n=n^{-1}\sum_i\delta_{X_i}$, independence and the basis bound give
$\E_f\left\|\theta_\ell\left(P_n\right)-\theta_\ell\left(f\right)\right\|_2^2\le B_0^2D_\ell/n$.
Uniformly convergent test expansions and Minkowski's inequality then yield
\[
\left(\E_f d_q\left(P_n,P_f\right)^2\right)^{1/2}
\le\sum_\ell\rho_\ell
 \left(\E_f\left\|\theta_\ell\left(P_n\right)-\theta_\ell\left(f\right)\right\|_2^2\right)^{1/2}
\le C_d/\sqrt n,
\]
since $\sum_\ell\rho_\ell\sqrt{D_\ell}<\infty$. Jensen gives the
first-moment bound.

For $f\in\mathcal F_\beta\left(R\right)$, the Fourier weights on $A_\ell$ are bounded below by a constant times $2^{2\beta\ell}$. Thus $\left\|P_\ell f\right\|_2\le C_\beta R2^{-\beta\ell}$. Summing the resulting geometric sequence gives, for every $1\le r\le\infty$,
\begin{equation}
\left\|\left(\rho_\ell\left\|P_\ell f\right\|_2\right)_{\ell\ge K}\right\|_{\ell_r}
\le C_{d,\beta}R2^{-\left(\beta+d\right)K}.
\label{eq:bias}
\end{equation}
The same bound applies to the difference of two members of the class with $R$ replaced by $2R$. This observation will allow one proper output rule to handle every fine index.

\subsection{Comparison with the earlier Sobolev class}
\label{app:albert}
We use the published convention in Section~2.1, equation~(5), of
\citet{AlbertEtAl2026}. Let
\[
 \varphi_1\left(t\right)=1,\qquad
 \varphi_{2j}\left(t\right)=\sqrt2\cos\left(2\pi jt\right),\qquad
 \varphi_{2j+1}\left(t\right)=\sqrt2\sin\left(2\pi jt\right),\quad j\ge1,
\]
and let $\varphi_{\boldsymbol j}\left(x\right)=\prod_{i=1}^d\varphi_{j_i}\left(x_i\right)$
for $\boldsymbol j\in\left\{1,2,\ldots\right\}^d$. For $\sigma>0$, define
\begin{equation*}
 \begin{gathered}
 \left\|u\right\|_{A,\sigma}^2
 =\sum_{\boldsymbol j\in\left\{1,2,\ldots\right\}^d}
       \left(\sum_{i=1}^d j_i^{2\sigma}\right)
       |\left\langle u,\varphi_{\boldsymbol j}\right\rangle|^2,\\
 \underline c_{d,\sigma}=\left(1+4\pi^2d\right)^{-\sigma},\qquad
 \overline C_{d,\sigma}=d\,3^{2\sigma}.
 \end{gathered}
\end{equation*}
In particular, the constant function one has squared $A$-norm $d$.
The identification of $L^2\left(\left[0,1\right]^d\right)$ with $L^2\left(\T^d\right)$ below uses
this same complete Fourier basis.

\begin{lemma}
\label{lem:albert-classes}
For every real, zero-mean $u\in L^2\left(\T^d\right)$ and every $\sigma>0$,
with infinite norms permitted,
\begin{equation}
 \underline c_{d,\sigma}\left\|u\right\|_{H^\sigma}^2
 \le\left\|u\right\|_{A,\sigma}^2
 \le \overline C_{d,\sigma}\left\|u\right\|_{H^\sigma}^2.
 \label{eq:albert-norm-comparison}
\end{equation}
For $\beta>0$ and $R_A>0$, write
\[
 \mathcal A_\beta\left(R_A\right)
 =\left\{f\in L^2\left(\T^d\right):f\ge0\ \text{a.e.},
       \int f=1,\ \left\|f\right\|_{A,\beta}\le R_A\right\}.
\]
If $R_A^2>d$, set
\begin{equation}
 R_-^2=1+\frac{R_A^2-d}{\overline C_{d,\beta}},\qquad
 R_+^2=1+\frac{R_A^2-d}{\underline c_{d,\beta}}.
 \label{eq:albert-radii}
\end{equation}
Then $R_-,R_+>1$ and
\begin{equation*}
 \mathcal F_\beta\left(R_-\right)
 \subseteq\mathcal A_\beta\left(R_A\right)
 \subseteq\mathcal F_\beta\left(R_+\right).
\end{equation*}
At $R_A^2=d$ the density class consists only of the uniform density;
at $R_A^2<d$ it is empty.

For a discriminator radius $R_D>0$, define
\[
 d_{A,R_D}\left(\mu,\nu\right)
 =\sup_{\left\|h\right\|_{A,d}\le R_D}
       \left|\int h\,d\left(\mu-\nu\right)\right|,
 \qquad \Lambda_d=\left(1+16\pi^2d\right)^{d/2}.
\]
For all probability measures $\mu,\nu$ on $\T^d$,
\begin{equation}
 \frac{R_D}{\sqrt{\overline C_{d,d}}\,\Lambda_d}\,d_2\left(\mu,\nu\right)
 \le d_{A,R_D}\left(\mu,\nu\right)
 \le\frac{R_D}{\sqrt{\underline c_{d,d}}}\,d_2\left(\mu,\nu\right).
 \label{eq:albert-loss-comparison}
\end{equation}
Here $d_2$ is the IPM with fine index $q=2$, not its square.
\end{lemma}

\begin{proof}
Group complex frequencies by $\left(|k_1|,\ldots,|k_d|\right)=\boldsymbol\nu$.
The corresponding real tensor indices have $j_i=1$ if $\nu_i=0$,
and $j_i\in\left\{2\nu_i,2\nu_i+1\right\}$ otherwise. The two orthonormal bases
span the same complexified subspace and preserve its coefficient energy.
Writing $M=\max\left\{1,|\boldsymbol\nu|_\infty\right\}$, their weights satisfy
\[
 M^{2\sigma}\le\sum_i j_i^{2\sigma}\le d3^{2\sigma}M^{2\sigma},\qquad
 M^{2\sigma}\le\left(1+4\pi^2|\boldsymbol\nu|_2^2\right)^\sigma
 \le\left(1+4\pi^2d\right)^\sigma M^{2\sigma}.
\]
Comparison on each subspace and summation prove
\eqref{eq:albert-norm-comparison}, including infinite norms.
For $f=1+u$, orthogonality gives
$\left\|f\right\|_{H^\beta}^2=1+\left\|u\right\|_{H^\beta}^2$ and
$\left\|f\right\|_{A,\beta}^2=d+\left\|u\right\|_{A,\beta}^2$.
Substitution of \eqref{eq:albert-radii} proves the class inclusions
and the assertions for $R_A^2\le d$.

For zero-mean $h$, let $B\left(h\right)=\left\|h\right\|_{\mathcal G_2}$.
On $A_\ell$, the $H^d$ weight lies between $2^{2d\ell}$ and
$\left(1+16\pi^2d\right)^d2^{2d\ell}$. Thus Plancherel gives
$B\left(h\right)\le\left\|h\right\|_{H^d}\le\Lambda_dB\left(h\right)$ and
\[
 \sqrt{\underline c_{d,d}}B\left(h\right)\le\left\|h\right\|_{A,d}
 \le\sqrt{\overline C_{d,d}}\Lambda_dB\left(h\right).
\]
Subtracting a test function's mean decreases both norms and preserves
its pairing with $\mu-\nu$. The resulting inclusions of test balls
prove \eqref{eq:albert-loss-comparison}. The block bounds above give
continuous representatives, so the comparison covers arbitrary
probability measures.
\end{proof}

For fixed $\beta>0$, $R_A>\sqrt d$, and $R_D>0$, the class inclusions and \eqref{eq:albert-loss-comparison} transfer Theorem~\ref{thm:critical} at $q=2$ to $\mathcal A_\beta\left(R_A\right)$: its minimax $p$th-moment risk under $d_{A,R_D}$ has order $R_D^p N^{-p/2}\left(\log N\right)^p$, $p=1,2$. Apply the lower bound on the inner ball $\mathcal F_\beta\left(R_-\right)$ and the upper bound on the outer ball $\mathcal F_\beta\left(R_+\right)$; the latter estimator takes values in the outer ball. Constants depend on the fixed class parameters, uniformly for $0<\varepsilon\le1$ and sufficiently large $N$.

\section{Proof of the upper bound}
\label{app:upper}

\subsection{Privacy and the spherical decoder}

By the basis bound and the ellipsoid definition \eqref{eq:ellipsoid},
\[
\left\|A^{-1}z\left(x\right)\right\|_2^2
=\frac1{B_0^2S_K}\sum_{\ell<K}\rho_\ell\left\|v_\ell\left(x\right)\right\|_2^2
\le1.
\]
Thus \eqref{eq:bit} is a valid kernel. Its ratio under any two inputs is at most $\left(1+\kappa\right)/\left(1-\kappa\right)=e^\varepsilon$, for every public $U_i$. Since the directions are independent of the inputs, the full noninteractive transcript satisfies \eqref{eq:privacy} for each user.

If $U$ is uniform on $\mathbb S^{m-1}$ and $t\ne0$, rotation invariance gives
\[
\E\left[U\sign\left\langle U,t\right\rangle\right]
=c_m\frac{t}{\left\|t\right\|_2}.
\]
Indeed the component parallel to $t$ has expectation $\E|U_1|$, and every orthogonal component has expectation zero by reflection. Applying this identity to \eqref{eq:bit} and \eqref{eq:decoder} yields
\[
\E\left[W_i\mid X_i=x\right]
=\frac{A}{\kappa c_{m_K}}
\E_U\left[U\kappa\left\|t\left(x\right)\right\|_2\sign\left\langle U,t\left(x\right)\right\rangle\right]
=At\left(x\right)=z\left(x\right).
\]
The identity is also valid at $t=0$. This proves unbiased decoding.

\subsection{Directional moments and block maxima}

For a scalar random variable $V$, use
$\left\|V\right\|_{\psi_2}=\inf\left\{s>0:\E\exp\left(V^2/s^2\right)\le2\right\}$. The spherical moment formula gives, for $j\ge1$,
\[
\E|U_1|^{2j}
=\frac{\left(2j-1\right)!!}{m\left(m+2\right)\cdots\left(m+2j-2\right)}
\le\frac{\left(2j\right)^j}{m^j}.
\]
The displayed formula follows, for example, by writing a standard Gaussian vector as its independent radius times a uniform direction and comparing even moments. It also gives the gamma function expression for $c_m$, with $c_m\asymp m^{-1/2}$. Expanding an exponential series, or summing the even-moment bounds, therefore yields
$\left\|\left\langle a,U\right\rangle\right\|_{\psi_2}\le C\left\|a\right\|_2/\sqrt m$.

For every deterministic vector $a\in\R^{m_K}$, the fact that $|B_i|=1$ gives
\[
|\left\langle a,W_i\right\rangle|
=\frac{|\left\langle A^\top a,U_i\right\rangle|}{\kappa c_{m_K}}.
\]
Hence the previous moment bound, centering, and $\kappa\asymp\varepsilon$ give
\begin{equation}
\left\|\left\langle a,W_i-\E_fW_i\right\rangle\right\|_{\psi_2}
\le C\frac{\left\|A^\top a\right\|_2}{\varepsilon}.
\label{eq:directional}
\end{equation}
This expectation includes the observation and both public and private randomization. The dependence of $B_i$ on $U_i$ does not affect the absolute-value identity.

A centered scalar variable with $\psi_2$ norm at most $s$ has moment-generating function bounded by $\exp\left(C\lambda^2s^2\right)$; this follows by symmetrization and expansion of the exponential using the even moments. Multiplying these bounds for independent users and optimizing $\lambda$ in Markov's inequality, \eqref{eq:directional} yields, for a unit vector $a$ in block $\ell$,
\[
\Prob_f\left(|\left\langle a,E_\ell\right\rangle|>t\right)
\le2\exp\left(-\frac{cNt^2}{S_K\rho_\ell}\right)
\le2\exp\left(-\frac{c'ND_\ell t^2}{K}\right).
\]
Here $S_K\asymp K$ and $D_\ell\rho_\ell\asymp1$.

Choose a $1/2$-net $\mathcal N_\ell$ of the Euclidean unit sphere in $\R^{D_\ell}$ with cardinality at most $5^{D_\ell}$. Such a net follows by taking a maximal separated set and comparing disjoint Euclidean balls. For every vector $v$, $\left\|v\right\|_2\le2\max_{a\in\mathcal N_\ell}|\left\langle a,v\right\rangle|$. The union bound therefore gives
\[
\Prob_f\left(\left\|E_\ell\right\|_2>A_0\sqrt{K/N}\,u\right)
\le2\exp\left(D_\ell\log5-cA_0^2D_\ell u^2\right)
\le2e^{-c_1D_\ell u^2}
\]
for $u\ge1$ after increasing the fixed constant $A_0$. Since $D_\ell\ge c_d2^{d\ell}$,
\[
\sum_{\ell\ge0}e^{-c_1D_\ell u^2}
\le e^{-c_2u^2}\sum_{\ell\ge0}e^{-c_2 2^{d\ell}}
\le C_de^{-c_2u^2},\qquad u\ge1,
\]
with a smaller $c_2>0$. This proves \eqref{eq:block-tail} without using independence between blocks. Integrating the tail, including the interval below $A_0\sqrt{K/N}$, proves
\begin{equation}
\sup_f\E_f\left(W_K^*\right)^p\le C_p\left(K/N\right)^{p/2},
\qquad p\in\left\{1,2\right\}.
\label{eq:max-moments}
\end{equation}
The supremum here may be taken over all input densities; smoothness is used only for the bias.

\subsection{Measurable density output and both moments}
\label{app:proper-output}

The subset $\mathcal F_\beta\left(R\right)$ of the separable space $L^2$ has a countable dense family $\left(g_m\right)$. Choose measurable representatives of its members. By \eqref{eq:l2-continuity}, $Y\mapsto d_1\left(Y,g_m\right)$ is continuous on the mass-one affine subspace of $L^2$. The countable infimum $I\left(Y\right)$ in \eqref{eq:selector} is Borel. The index $m\left(Y\right)$ is Borel because the event $m\left(Y\right)=m$ is the event that the $m$th inequality holds and all earlier ones fail. The set of eligible indices is nonempty since $\tau>0$.

Density of the family and \eqref{eq:l2-continuity} imply $I\left(Y\right)\le d_1\left(Y,f\right)$ for every $f\in\mathcal F_\beta\left(R\right)$. The triangle inequality gives
\begin{equation}
d_1\left(\widehat f,f\right)\le2d_1\left(Y,f\right)+\tau
\le2W_K^*+C R2^{-\left(\beta+d\right)K}+\tau,
\label{eq:proper-max}
\end{equation}
where \eqref{eq:block-duality}, \eqref{eq:decoder}, and \eqref{eq:bias} give the last bound. The positive tolerance ensures that the countable family contains an eligible index, yielding a measurable approximate minimum-distance estimate.

Now both $\widehat f$ and $f$ lie in $\mathcal F_\beta\left(R\right)$. Split their difference at $K$, apply \eqref{eq:aggregation} to the retained blocks, and apply \eqref{eq:bias} to the remainder. For $r=q'$, this gives
\[
\begin{split}
d_q\left(\widehat f,f\right)
&\le K^{1/r}d_1\left(\widehat f,f\right)+CR2^{-\left(\beta+d\right)K}\\
&\le 2K^{1/r}W_K^*
+C_RK^{1/r}\left(2^{-\left(\beta+d\right)K}+\tau\right),
\end{split}
\]
using \eqref{eq:proper-max} and $K\ge1$. The cutoff \eqref{eq:cutoff} gives $2^{-\left(\beta+d\right)K}\le N^{-1/2}=\tau$. Taking the first or second moment and using \eqref{eq:max-moments} yields
\[
\sup_{f\in\mathcal F_\beta\left(R\right)}\E_f d_q\left(\widehat f,f\right)^p
\le C_{p,d,\beta,R,q}\frac{K^{p\left(1/2+1/r\right)}}{N^{p/2}},
\qquad p=1,2.
\]
Because $K\asymp\log N$, this proves the upper half of Theorem~\ref{thm:critical}. In particular, the squared-risk upper uses the actual second-moment bound in \eqref{eq:max-moments}.

\subsection{Known radius and unknown smoothness}
\label{app:smoothness-free-proof}

\begin{proof}[Proof of Corollary~\ref{cor:smoothness-free}]
Write $Y_0$ for the Fourier estimate in \eqref{eq:decoder} with cutoff $K_0$. The privacy, unbiasedness, and measurability arguments above apply with this cutoff and the separable class $\mathcal F_0\left(R\right)$. Equation~\eqref{eq:max-moments} gives $\E_f\left(W_{K_0}^*\right)^p\le C_p\left(K_0/N\right)^{p/2}$ for $p=1,2$, uniformly over input densities. It remains to bound the omitted blocks using only the $L^2$ radius.

For every mass-zero $u\in L^2$ and $1\le r\le\infty$,
\begin{equation}
 \left\|\left(\rho_\ell\left\|P_\ell u\right\|_2\right)_{\ell\ge K_0}\right\|_{\ell_r}
 \le\frac{\left\|u\right\|_2\,2^{-dK_0}}{1-2^{-d}}.
 \label{eq:smoothness-free-tail}
\end{equation}
For $u=\widehat f^{\left(0\right)}-f$ this is at most $4RN^{-1/2}$.
The retained blocks of $Y_0-f$ have weighted norms $\left\|E_\ell\right\|_2$;
the omitted blocks give
$d_1\left(Y_0,f\right)\le W_{K_0}^*+R2^{-dK_0}$.
Density of the candidate family and the selector inequality imply
$d_1\left(\widehat f^{\left(0\right)},f\right)\le2W_{K_0}^*+2R2^{-dK_0}+N^{-1/2}$.
Splitting \eqref{eq:block-duality} at $K_0$, applying
\eqref{eq:aggregation} to the retained blocks and
\eqref{eq:smoothness-free-tail} to the rest, yields
\[
 d_q\left(\widehat f^{\left(0\right)},f\right)
 \le2K_0^{1/r}W_{K_0}^*+\left(6R+1\right)K_0^{1/r}N^{-1/2},\qquad r=q'.
\]
Taking first and second moments using
\eqref{eq:max-moments} proves
\eqref{eq:smoothness-free-risk}. Since $K_0\asymp_d\log N$ and
$0\le1/r\le1$, its constant can be chosen uniformly in $q$.

Every $\mathcal F_\beta\left(R\right)$ is contained in $\mathcal F_0\left(R\right)$.
For $\beta>0$, Theorem~\ref{thm:critical} supplies the matching lower
bound; for $\beta=0$, restrict to $\mathcal F_1\left(R\right)\subset\mathcal F_0\left(R\right)$.
The same protocol and output classes are allowed in each lower bound.
Finally, $m_K=\left(2^{K+1}-1\right)^d-1\asymp_d2^{dK}$ gives the release dimensions.
\end{proof}

\subsection{Finite-dimensional output in the \texorpdfstring{$L^2$}{L2} ball}

\begin{proposition}
\label{prop:f0-finite}
Fix $d\ge1$ and a known $R>1$. For $0<\varepsilon\le1$ and
$N=n\varepsilon^2\ge2$, use the release
\eqref{eq:ellipsoid}--\eqref{eq:decoder} with
\[
 K_0=\left\lceil\frac{\log_2N}{2d}\right\rceil,\qquad
 L=\left\lceil\sqrt N\right\rceil,\qquad
 \delta=N^{-1/2},\qquad \eta=\delta/R^2.
\]
Partition $\T^d$ into the $L^d$ half-open cubes $\mathcal P_L$ of
side length $1/L$, and define
\[
 \mathcal V_{L,R}
 =\left\{\boldsymbol v\in\R^{L^d}:
 v_C\ge0,\quad L^{-d}\sum_{C\in\mathcal P_L}v_C=1,\quad
 L^{-d}\sum_{C\in\mathcal P_L}v_C^2\le R^2\right\}.
\]
Write $h_{\boldsymbol v}=\sum_{C\in\mathcal P_L}v_C\mathbf1_C$.
For the decoded vector $\overline W$, minimize the finite-dimensional
convex objective
\begin{equation}
 \begin{split}
 \Phi_L(\boldsymbol v;\overline W)
 &=\max_{0\le\ell<K_0}
       \left\|\rho_\ell\theta_\ell(h_{\boldsymbol v})
                         -\overline W_\ell\right\|_2,\\
 J_L(\boldsymbol v;\overline W)
 &=\Phi_L(\boldsymbol v;\overline W)
                         +\eta L^{-d}\sum_Cv_C^2
 \end{split}
 \label{eq:f0-finite-objective}
\end{equation}
over $\mathcal V_{L,R}$. The exact minimizer is unique and depends
continuously on $\overline W$. More generally, let
$\widetilde{\boldsymbol v}$ be any measurable feasible output satisfying
\begin{equation}
 J_L(\widetilde{\boldsymbol v};\overline W)
 \le\min_{\boldsymbol v\in\mathcal V_{L,R}}
                   J_L(\boldsymbol v;\overline W)+\delta
 \label{eq:f0-finite-tolerance}
\end{equation}
almost surely, and set $\widehat f_L=h_{\widetilde{\boldsymbol v}}$.
Then $\widehat f_L\in\mathcal F_0(R)$, and a constant depending only on
$d,R$ satisfies, simultaneously for $1\le q\le\infty$ and $p=1,2$,
\begin{equation}
 \sup_{f\in\mathcal F_0(R)}
       \E_f d_q(\widehat f_L,f)^p
 \le C_{d,R}N^{-p/2}(\log N)^{p(3/2-1/q)}.
 \label{eq:f0-finite-risk}
\end{equation}
Thus the same finite-dimensional postprocessing attains the minimax
rates on every fixed $\mathcal F_\beta(R)$, $\beta\ge0$.
Its output is proper in $\mathcal F_0(R)$; no $H^\beta$ constraint
is imposed on the histogram. The program has $L^d\asymp N^{d/2}$
histogram variables, while the release dimension remains
$m_{K_0}\asymp\sqrt N$.
\end{proposition}

\begin{proof}
The feasible set is nonempty, compact, and convex: the vector with all
entries equal to one is feasible, and the quadratic constraint bounds
every coordinate. Its constraints give nonnegativity, integral one,
and $\|h_{\boldsymbol v}\|_2\le R$ exactly. The coefficient maps
$\boldsymbol v\mapsto\theta_\ell(h_{\boldsymbol v})$ are linear, with
entries $\int_C\psi_{\ell j}$. Consequently $\Phi_L$ is continuous
and convex. The positive quadratic term makes $J_L$ strictly convex,
so its minimizer exists and is unique. If decoded vectors converge,
the corresponding objectives converge uniformly on the compact
feasible set, since
\[
 |J_L(\boldsymbol v;w)-J_L(\boldsymbol v;w')|
 \le\max_{\ell<K_0}\|w_\ell-w'_\ell\|_2.
\]
Every convergent subsequence of minimizers therefore minimizes the
limiting objective. Uniqueness proves continuity of the minimizer,
and hence supplies a measurable feasible output satisfying
\eqref{eq:f0-finite-tolerance}. An epigraph formulation of
\eqref{eq:f0-finite-objective} uses one Euclidean-norm constraint for
each retained block, together with the displayed linear and
quadratic feasibility constraints.

Let $\Pi_L$ denote cell averaging and set
$\boldsymbol v_f=(L^d\int_C f)_{C\in\mathcal P_L}$, so that
$h_{\boldsymbol v_f}=\Pi_L f$. For $f\in\mathcal F_0(R)$, this function
is nonnegative and has integral one. Jensen's inequality on each
cell gives $\|\Pi_L f\|_2\le\|f\|_2\le R$, so it is a feasible
comparison density.

We next bound its Fourier approximation error without smoothness
assumptions on $f$. Cell averaging is a self-adjoint orthogonal
projection in $L^2$. For a continuously differentiable function
$u$ on an interval $I$ of length $1/L$, writing $u_I$ for its
average gives
\[
 \int_I|u-u_I|^2
 =\frac1{2|I|}\int_{I^2}|u(x)-u(y)|^2\,dx\,dy
 \le L^{-2}\int_I|u'|^2.
\]
The last inequality follows from the fundamental theorem of
calculus and Cauchy--Schwarz. Applying this estimate successively
in the $d$ coordinates yields
$\|h-\Pi_L h\|_2\le C_dL^{-1}\|\nabla h\|_2$.
For $h$ in the range of $P_\ell$, Plancherel gives
$\|\nabla h\|_2\le C_d2^\ell\|h\|_2$. Therefore duality gives
\begin{equation}
 \begin{split}
 \rho_\ell\|P_\ell(\Pi_L f-f)\|_2
 &=\rho_\ell
   \sup_{\substack{h=P_\ell h\\\|h\|_2\le1}}
          |\langle f,\Pi_L h-h\rangle|\\
 &\le C_d R L^{-1}2^{(1-d)\ell}
 \le C_dR/L.
 \end{split}
 \label{eq:f0-grid-approximation}
\end{equation}

Write $W_{K_0}^*=\max_{\ell<K_0}
\|\overline W_\ell-\rho_\ell\theta_\ell(f)\|_2$.
The feasible comparison density and
\eqref{eq:f0-grid-approximation} imply
$\Phi_L(\boldsymbol v_f;\overline W)\le W_{K_0}^*+C_dR/L$.
The regularization term is at most $\eta R^2=\delta$ on the
feasible set. Thus \eqref{eq:f0-finite-tolerance} gives
\[
 \max_{\ell<K_0}\rho_\ell
             \|P_\ell(\widehat f_L-f)\|_2
 \le 2W_{K_0}^*+C_dR/L+2\delta.
\]
Both densities have $L^2$ norm at most $R$. Splitting
\eqref{eq:block-duality} at $K_0$ and using
\eqref{eq:aggregation} and \eqref{eq:smoothness-free-tail} therefore
gives, for $r=q'$,
\[
 d_q(\widehat f_L,f)
 \le K_0^{1/r}\left(2W_{K_0}^*+C_dR/L+2\delta\right)
       +\frac{2R\,2^{-dK_0}}{1-2^{-d}}.
\]
Now $L^{-1}\le N^{-1/2}$ and $2^{-dK_0}\le N^{-1/2}$.
Taking the first and second moments using \eqref{eq:max-moments}
proves \eqref{eq:f0-finite-risk}, uniformly in $q$.
The matching lower bounds follow from Theorem~\ref{thm:critical}
for $\beta>0$ and from
$\mathcal F_1(R)\subseteq\mathcal F_0(R)$ for $\beta=0$.
The construction is postprocessing of the existing release and
therefore preserves its privacy guarantee.
\end{proof}
\section{Fisher information and a coordinate van Trees inequality}
\label{app:information}

We give the information argument for the public blackboard experiment defined above. The contraction is the bounded score-covariance case of \citet{BarnesChenOzgur2020}, including their blackboard extension. The proof below also records the regularity needed for the coordinate form of the classical van Trees inequality~\citep{GillLevit1995}.

\begin{lemma}
\label{lem:fisher}
Let $\psi_1,\ldots,\psi_m$ be bounded, orthonormal in $L^2\left(\mu\right)$ with $\int\psi_j\,d\mu=0$, where $\left(\mathcal X,\mu\right)$ is a standard Borel probability space, such as $\T^d$ or $\left[0,1\right]$ with Lebesgue measure, and let
\[
 f_\theta=1+\sum_{j=1}^m\theta_j\psi_j,
 \qquad \theta\in\Theta,
\]
where $\Theta$ is a compact rectangle on which $1/2\le f_\theta\le3/2$. For $n$ independent observations with density $f_\theta$, let $\mathsf Z$ be the complete transcript of any protocol in $\mathcal Q_\varepsilon^{\mathrm{BB}}\left(n\right)$. Its Fisher information matrix satisfies
\begin{equation}
 \tr I_{\mathsf Z}\left(\theta\right)
 \le 2n\left(e^\varepsilon-1\right)^2=2N_\varepsilon.
 \label{eq:private-information}
\end{equation}
The same bound holds if the protocol also involves independent users whose laws do not depend on $\theta$: only the $n$ users carrying information about $\theta$ enter the right-hand side.
\end{lemma}

\begin{proof}
Write $P_\theta$ for the law of one informative observation. Its score is
\[
 s_\theta\left(x\right)=\frac{\left(\psi_1\left(x\right),\ldots,\psi_m\left(x\right)\right)^{\mathsf T}}
 {f_\theta\left(x\right)}.
\]
The mean of this score is zero. Moreover, for every $v\in\R^m$,
\[
 \E_\theta\left\langle v,s_\theta\left(X\right)\right\rangle^2
 =\int_{\mathcal X}\frac{\left(\sum_jv_j\psi_j\left(x\right)\right)^2}{f_\theta\left(x\right)}\,d\mu\left(x\right)
 \le2\left\|v\right\|_2^2.
\]

Let $\nu$ be the transcript law when the informative inputs have density one relative to $\mu$, and write $\E_0$ for expectation in this reference experiment. Any parameter-independent users retain their fixed laws. The parameter-independent protocol kernel gives the likelihood
\[
 p_\theta\left(z\right)=\frac{dP_\theta^{\mathsf Z}}{d\nu}\left(z\right)
 =\E_0\!\left[\prod_{i=1}^nf_\theta\left(X_i\right)\,\middle|\,\mathsf Z=z\right].
\]
For fixed $n,m$, this is a polynomial in $\theta$ with bounded
measurable coefficients. On one full-measure set,
$2^{-n}\le p_\theta\le\left(3/2\right)^n$ throughout the compact rectangle.
Differentiation under the integral is valid, giving mean-zero scores and
\begin{equation}
 \nabla_\theta\log p_\theta\left(\mathsf Z\right)
 =\sum_{i=1}^n \bar s_i\left(\mathsf Z\right),
 \qquad
 \bar s_i\left(\mathsf Z\right)=\E_\theta\left[s_\theta\left(X_i\right)\mid\mathsf Z\right].
 \label{eq:transcript-score}
\end{equation}
The same argument gives finite second moments of every coordinate score. Positivity on a neighborhood of $\Theta$, if needed for derivatives at its boundary, follows by continuity and compactness of the finite-dimensional family.

Conditional on the public seed, the initial input law is a product.
A message multiplies only its speaker's posterior factor, so induction
over the finite transcript preserves conditional independence. Repeated
visits update the same factor; public scheduling and messages independent
of private inputs preserve the induction. This conditional-distribution argument
on standard Borel spaces requires no common message density.

For a fixed $\theta$, integrate the privacy inequality~\eqref{eq:privacy} over all inputs except that of user $i$. This gives an $\varepsilon$-private marginal channel
\[
 K_{i,\theta}\left(B\mid x\right)
 =\int Q\left(B\mid x,x_{-i}\right)\,dP_\theta^{n-1}\left(x_{-i}\right).
\]
Any additional parameter-independent users are integrated out as well. Averaging the pairwise privacy inequalities over a second input $x'\sim P_\theta$ gives, as inequalities of measures,
\[
 e^{-\varepsilon}P_\theta^{\mathsf Z}
 \le K_{i,\theta}\left(\cdot\mid x\right)
 \le e^\varepsilon P_\theta^{\mathsf Z}.
\]
Consequently the posterior has a density
\[
 \rho_{i,\theta}\left(z,x\right)
 =\frac{dP_\theta^{X_i\mid\mathsf Z=z}}{dP_\theta}\left(x\right),
 \qquad
 e^{-\varepsilon}\le\rho_{i,\theta}\left(z,x\right)\le e^\varepsilon
\]
for $P_\theta^{\mathsf Z}\otimes P_\theta$-almost every $\left(z,x\right)$. This conclusion applies privacy to the joint transcript, including the public seed, before disintegration.

For $\bar s_i\left(z\right)\ne0$, take $v=\bar s_i\left(z\right)/\left\|\bar s_i\left(z\right)\right\|_2$. The zero mean of the unprivatized score and Cauchy-Schwarz give
\begin{align*}
 \left\|\bar s_i\left(z\right)\right\|_2^2
 &=\left(\int\left\langle v,s_\theta\left(x\right)\right\rangle
       \left(\rho_{i,\theta}\left(z,x\right)-1\right)\,dP_\theta\left(x\right)\right)^2\\
 &\le2\int\left(\rho_{i,\theta}\left(z,x\right)-1\right)^2\,dP_\theta\left(x\right)
 \le2\left(e^\varepsilon-1\right)^2.
\end{align*}
The same inequality is immediate when $\bar s_i\left(z\right)=0$. Conditional independence gives, for $i\ne k$,
\[
 \E_\theta\left[\bar s_i\bar s_k^{\mathsf T}\right]
 =\E_\theta\!\left[
   \E_\theta\left[s_\theta\left(X_i\right)s_\theta\left(X_k\right)^{\mathsf T}\mid\mathsf Z\right]
   \right]
 =\E_\theta\left[s_\theta\left(X_i\right)s_\theta\left(X_k\right)^{\mathsf T}\right]=0.
\]
Combining this identity with~\eqref{eq:transcript-score} proves~\eqref{eq:private-information}. Parameter-independent users have zero score, which proves the final assertion.
\end{proof}

\begin{lemma}
\label{lem:van-trees}
For an experiment with the likelihood regularity established in Lemma~\ref{lem:fisher}, put $\Theta=\prod_{j=1}^m\left[-a_j,a_j\right]$ with $a_j>0$, and equip $\Theta$ with the product prior $\Pi$ whose coordinate densities are
\[
 \pi_{a_j}\left(t\right)=\frac1{a_j}
 \cos^2\!\left(\frac{\pi t}{2a_j}\right)
 \mathbf1_{\left\{|t|\le a_j\right\}}.
\]
Let $\E_\Pi$ denote expectation over both the prior and the experiment, and define
$b_j=\int_\Theta\left[I_{\mathsf Z}\left(\theta\right)\right]_{jj}\,d\Pi\left(\theta\right)$.
For every bounded measurable coordinate estimate $\widehat\theta_j$,
\begin{equation}
 \E_\Pi\left(\widehat\theta_j-\theta_j\right)^2
 \ge\frac1{b_j+\pi^2/a_j^2}.
 \label{eq:coordinate-van-trees}
\end{equation}
If $\widehat\theta_j\in\left[-a_j,a_j\right]$, then also
\begin{equation}
 \E_\Pi|\widehat\theta_j-\theta_j|
 \ge\frac1{2a_j\left(b_j+\pi^2/a_j^2\right)}.
 \label{eq:coordinate-first-moment}
\end{equation}
\end{lemma}

\begin{proof}
The prior integrates to one, vanishes at its endpoints, and has Fisher
information $\pi^2/a_j^2$. Write $\pi_\Theta(\theta)=\prod_{k=1}^m\pi_{a_k}(\theta_k)$ for the joint prior density. For the joint score
$V_j=\partial_j\log p_\theta\left(\mathsf Z\right)+\partial_j\log\pi_\Theta\left(\theta\right)$,
integration by parts gives
\[
\E_\Pi\left[\left(\widehat\theta_j-\theta_j\right)V_j\right]=1,\qquad
\E_\Pi V_j^2=b_j+\pi^2/a_j^2.
\]
The vanishing prior, bounded estimate and established likelihood regularity
justify integration by parts. The cross term in the second identity is
zero because the transcript score has conditional mean zero.
Cauchy-Schwarz proves \eqref{eq:coordinate-van-trees};
$\left(\widehat\theta_j-\theta_j\right)^2\le2a_j|\widehat\theta_j-\theta_j|$
proves \eqref{eq:coordinate-first-moment}. Independent estimator
randomness is parameter-independent postprocessing.
\end{proof}

\section{The multiscale lower bound}
\label{app:lower}

We prove the lower bounds in Theorem~\ref{thm:critical}. Fix $c_0>0$ and choose
\[
 0<\gamma<\min\left\{\frac1{3d},\frac1{2\beta+2d}\right\},
 \qquad K=\left\lfloor\gamma\log_2N_\varepsilon\right\rfloor.
\]
For sufficiently large $N$, this is a positive integer. In each shell $A_\ell$, $0\le\ell<K$, retain all $D_\ell$ real basis functions $\psi_{\ell j}$. Set
\begin{equation}
 a_\ell^2=\frac{c_0KD_\ell}{N_\varepsilon},
 \qquad
 f_\theta=1+\sum_{\ell=0}^{K-1}\sum_{j=1}^{D_\ell}
                 \theta_{\ell j}\psi_{\ell j},
 \qquad |\theta_{\ell j}|\le a_\ell.
 \label{eq:lower-prior-widths}
\end{equation}
Give these coordinates independent priors $\pi_{a_\ell}$ from Lemma~\ref{lem:van-trees}.

All basis functions have mean zero, so every $f_\theta$ integrates to one. Their uniform boundedness and $D_\ell\asymp2^{d\ell}$ imply
\begin{align*}
 \left\|f_\theta-1\right\|_\infty
 &\le C\sum_{\ell<K}D_\ell a_\ell
 \le C\sqrt{\frac K{N_\varepsilon}}\,2^{3dK/2},\\
 \left\|f_\theta\right\|_{H^\beta}^2-1
 &\le C\sum_{\ell<K}2^{2\beta\ell}D_\ell a_\ell^2
 \le C\frac K{N_\varepsilon}\,2^{\left(2\beta+2d\right)K}.
\end{align*}
Both bounds tend to zero by the strict inequalities in the choice of $\gamma$. Since $R>1$, the entire parameter rectangle therefore belongs to $\mathcal F_\beta\left(R\right)$ and satisfies $1/2\le f_\theta\le3/2$ once $N$ is sufficiently large. The upper and lower constructions may use different constant multiples of $\log N$ for their respective cutoffs.

Fix any permitted protocol and any measurable estimate $\widehat\mu$. For probability-measure-valued estimates, each selected coefficient
$\widetilde\theta_{\ell j}=\int\psi_{\ell j}\,d\widehat\mu$
is well defined and measurable. The same reasoning applies to any other output for which these finite Fourier pairings and the loss are defined. Clip $\widetilde\theta_{\ell j}$ to $\left[-a_\ell,a_\ell\right]$ and denote the resulting estimate by $\widehat\theta_{\ell j}$. Put
\[
 e_{\ell j}=\widehat\theta_{\ell j}-\theta_{\ell j},
 \qquad
 t_\ell=\rho_\ell
       \left(\sum_{j=1}^{D_\ell}e_{\ell j}^2\right)^{1/2}.
\]
Clipping decreases every coordinate error in absolute value. The finite-block loss inequality~\eqref{eq:finite-duality} consequently gives
\[
 d_q\left(\widehat\mu,f_\theta\right)
 \ge\left\|\left(t_\ell\right)_{\ell<K}\right\|_{\ell_r},
 \qquad \frac1q+\frac1r=1.
\]

Write
$b_{\ell j}=\int\left[I_{\mathsf Z}\left(\theta\right)\right]_{\ell j,\ell j}\,d\Pi\left(\theta\right)$.
The information bound~\eqref{eq:private-information} and the widths in~\eqref{eq:lower-prior-widths} give the two exact budget bounds
\begin{equation}
 \sum_{\ell<K}\sum_{j=1}^{D_\ell}b_{\ell j}\le2N_\varepsilon,
 \qquad
 \sum_{\ell<K}D_\ell\frac{\pi^2}{a_\ell^2}
 =\frac{\pi^2N_\varepsilon}{c_0}.
 \label{eq:lower-information-budgets}
\end{equation}
Applying~\eqref{eq:coordinate-first-moment} and using
$\left\|e_\ell\right\|_2\ge D_\ell^{-1/2}\left\|e_\ell\right\|_1$, we obtain
\[
 \E_\Pi\sum_{\ell<K}t_\ell
 \ge\sum_{\ell<K}\sum_{j=1}^{D_\ell}
 \frac{u_\ell}{b_{\ell j}+\pi^2/a_\ell^2},
 \qquad
 u_\ell=\frac{\rho_\ell}{2a_\ell\sqrt{D_\ell}}.
\]
Weighted Cauchy-Schwarz now yields
\begin{equation}
 \E_\Pi\sum_{\ell<K}t_\ell
 \ge
 \frac{\left(\sum_{\ell<K}D_\ell\sqrt{u_\ell}\right)^2}
 {\sum_{\ell<K}\sum_{j=1}^{D_\ell}b_{\ell j}
       +\sum_{\ell<K}D_\ell\pi^2/a_\ell^2}.
 \label{eq:lower-weighted-cauchy}
\end{equation}
The geometric cancellation at the critical loss can be seen explicitly:
\[
 D_\ell\sqrt{u_\ell}
 =\frac{1}{\sqrt2}\left(\frac{N_\varepsilon}{c_0K}\right)^{1/4}
   \sqrt{D_\ell\rho_\ell}
 \asymp\left(\frac{N_\varepsilon}{K}\right)^{1/4}.
\]
Thus the numerator in~\eqref{eq:lower-weighted-cauchy} is bounded below by a constant times $N_\varepsilon^{1/2}K^{3/2}$. Its denominator is at most $\left(2+\pi^2/c_0\right)N_\varepsilon$ by~\eqref{eq:lower-information-budgets}. It follows that
\[
 \E_\Pi\sum_{\ell<K}t_\ell
 \ge c\frac{K^{3/2}}{\sqrt{N_\varepsilon}}.
\]

For every nonnegative vector of length $K$ and every $1\le r\le\infty$,
\[
 \left\|t\right\|_{\ell_r}\ge K^{1/r-1}\left\|t\right\|_{\ell_1}.
\]
The preceding bounds therefore imply, uniformly over all allowed protocols and estimates,
\[
 \E_\Pi d_q\left(\widehat\mu,f_\theta\right)
 \ge c\frac{K^{1/2+1/r}}{\sqrt{N_\varepsilon}}.
\]
In particular, $r=\infty$ uses the factor $K^{-1}$ and gives the maximum-block lower bound $c\sqrt{K/N_\varepsilon}$. No additional logarithm is introduced at this endpoint. For the squared loss, Jensen's inequality on the joint prior-and-transcript law gives
\[
 \E_\Pi d_q\left(\widehat\mu,f_\theta\right)^2
 \ge\left(\E_\Pi d_q\left(\widehat\mu,f_\theta\right)\right)^2
 \ge c\frac{K^{1+2/r}}{N_\varepsilon}.
\]
The supremum risk over $\mathcal F_\beta\left(R\right)$ dominates each of these Bayes risks, since the entire prior is supported in that class. Finally,
$N\le N_\varepsilon\le\left(e-1\right)^2N$ for $0<\varepsilon\le1$,
$K\asymp\log N$, and $1/2+1/r=3/2-1/q$.
Taking the infimum over protocols and estimators proves both lower bounds in Theorem~\ref{thm:critical}.

\section{Proofs for cumulative distributions and transport}\label{app:cdf}

\subsection{A hierarchical CDF estimator}
\label{app:cdf-hierarchy}

We give a complete population version of the hierarchical construction of
\citet[Sections~3.2, 4.4, and 4.7]{CormodeKulkarniSrivastava2019}.
The frequency oracle is optimized unary encoding, and the levels use disjoint,
deterministically assigned groups of users. Throughout this appendix,
unsubscripted $L^2$ norms refer to Lebesgue measure on $\left[0,1\right]$.

Let $h=\left\lceil\left(\log_2 N\right)/2\right\rceil$ and $D=2^h$. Partition the users into
$h$ groups $\mathcal I_1,\ldots,\mathcal I_h$ whose sizes $n_\ell$ differ
by at most one. For sufficiently large $N$, uniformly over
$0<\varepsilon\le1$, we have $n_\ell\ge n/\left(2h\right)$, since $n\ge N$.
At level $\ell$, partition $\left[0,1\right]$ into $2^\ell$ dyadic cells
$I_{\ell,k}=\left[k2^{-\ell},\left(k+1\right)2^{-\ell}\right)$, closing the last cell at $1$.
For $i\in\mathcal I_\ell$, let $j$ be the cell containing $X_i$.
The user sends a vector of conditionally independent bits
$V_i=\left(V_{i,k}:0\le k<2^\ell\right)$ with
\[
 p_0=\frac12,\qquad q_0=\frac1{1+e^\varepsilon},\qquad
 \Prob\left(V_{i,k}=1\mid X_i\right)=
 \begin{cases}p_0,&k=j,\\q_0,&k\ne j.\end{cases}
\]
Changing $j$ to another cell changes only two Bernoulli factors. The largest
ratio of their joint probabilities is
$p_0\left(1-q_0\right)/\left[q_0\left(1-p_0\right)\right]=e^\varepsilon$; hence the full vector satisfies
pure $\varepsilon$-LDP. Every user sends one such vector, so the complete
protocol has total budget $\varepsilon$ per user. Each report at level $\ell$ contains $2^\ell$ bits, at most $D$.

Decode each cell probability by
\[
 \widehat p_{\ell,k}
 =\frac1{n_\ell}\sum_{i\in\mathcal I_\ell}
       \frac{V_{i,k}-q_0}{p_0-q_0}.
\]
Since $\E V_{i,k}=q_0+\left(p_0-q_0\right)P_f\left(I_{\ell,k}\right)$, this estimate is
unbiased. Its variance satisfies
\[
 \Var\left(\widehat p_{\ell,k}\right)
 \le\frac1{4n_\ell\left(p_0-q_0\right)^2}
 \le\frac{C}{n_\ell\varepsilon^2}
 \le\frac{Ch}{N}.
\]
This bound includes both sampling variation and privatization noise.

For $0<j<D$, decompose $\left[0,j/D\right)$ into its maximal dyadic subintervals.
There is at most one selected cell at each level. Let $\widehat F_j$ be
the sum of the corresponding cell estimates, and set
$\widehat F_0=0$, $\widehat F_D=1$. Different levels use independent
groups of observations and independent private coins. Therefore
\[
 \E_f\widehat F_j=F_f\left(j/D\right),\qquad
 \Var_f\left(\widehat F_j\right)\le\frac{Ch^2}{N}.
\]
No independence between estimates of different cells at the same level is
needed. Define $\widetilde F\left(x\right)=\widehat F_j$ on
$\left[j/D,\left(j+1\right)/D\right)$ and $\widetilde F\left(1\right)=1$.
When $f\le b$, the discretization bias on this interval is at most $b/D$.
Integrating the pointwise bias-variance decomposition gives
\begin{equation}\label{eq:cdf-hierarchy-bound}
 \E_f\left\|\widetilde F-F_f\right\|_2^2
 \le \frac{Ch^2}{N}+\frac{b^2}{D^2}
 \le C_b\frac{\left(\log N\right)^2}{N}.
\end{equation}
Only individual prefix mean squared errors enter this calculation.

\subsection{Projection and inversion}
\label{app:cdf-projection}

The convex set $\mathcal C_{a,b}$ defined by \eqref{eq:cdf-increments}
contains the identity and is compact
in $C\left[0,1\right]$ by Arzel\`a-Ascoli. Its continuous injection into $L^2$
has compact, hence closed image. Hilbert projection is unique,
nonexpansive and measurable, and for $F_f\in\mathcal C_{a,b}$,
\[
\left\|\widehat F-F_f\right\|_2^2+\left\|\widetilde F-\widehat F\right\|_2^2
\le\left\|\widetilde F-F_f\right\|_2^2.
\]
This proves \eqref{eq:cdf-projection-risk}; combined with
\eqref{eq:cdf-hierarchy-bound}, it proves the CDF upper bound.
The increment and endpoint conditions give an absolutely continuous CDF
with density in $\left[a,b\right]$.

The compact injection has continuous inverse onto its image, so selecting
the continuous representative is continuous from $L^2$ to $C\left[0,1\right]$.
The slope bound gives
$\left\|G^{-1}-H^{-1}\right\|_\infty\le a^{-1}\left\|G-H\right\|_\infty$.
Thus inversion and composition yield measurable map estimates.

Let $\mathcal C_{a,b}^{D}$ consist of the linear interpolants through
$\left(j/D,g_j\right)$, where $g_0=0$, $g_D=1$, and
$a/D\le g_{j+1}-g_j\le b/D$. This is a nonempty compact convex
subset of $\mathcal C_{a,b}$. Its exact $L^2$ projection
$\widehat F_D$ of $\widetilde F$ is a finite-dimensional quadratic
minimizer. The interpolant $F^D$ of the true CDF lies in this subset
and satisfies $\left\|F^D-F_f\right\|_\infty\le b/D$. Nonexpansiveness gives
\[
 \left\|\widehat F_D-F_f\right\|_2
 \le\left\|\widehat F_D-F^D\right\|_2+\left\|F^D-F_f\right\|_2
 \le\left\|\widetilde F-F_f\right\|_2+2b/D.
\]
Thus this finite-dimensional alternative also attains the rate in
Theorem~\ref{thm:cdf} and satisfies the same slope constraints.

\begin{proposition}
\label{prop:cdf-projection-gap}
Fix $0<a<1<b$, an integer $D\ge2$, and
$\widetilde F\in L^2\left[0,1\right]$. Let
\[
 \mathcal K_D=\left\{\boldsymbol g=\left(g_0,\ldots,g_D\right):
     g_0=0,\ g_D=1,
     a/D\le g_{j+1}-g_j\le b/D\ \left(0\le j<D\right)\right\},
\]
and let $G_{\boldsymbol g}$ be its linear interpolant. Define the
\emph{unhalved} squared objective and its linearization gap by
\begin{equation}
 \begin{split}
 J_D\left(\boldsymbol g\right)&=\left\|G_{\boldsymbol g}-\widetilde F\right\|_2^2,\\
 \operatorname{gap}_D\left(\boldsymbol g\right)
 &=\max_{\boldsymbol v\in\mathcal K_D}
       \sum_{i=1}^{D-1}\partial_iJ_D\left(\boldsymbol g\right)\left(g_i-v_i\right).
 \end{split}
 \label{eq:cdf-optimization-gap}
\end{equation}
Derivatives here are with respect to the free interior coordinates.
Let $\boldsymbol g_*$ minimize $J_D$ on $\mathcal K_D$, so that
$G_{\boldsymbol g_*}=\widehat F_D$. For every feasible
$\boldsymbol g$,
\begin{equation}
 \left\|G_{\boldsymbol g}-\widehat F_D\right\|_2^2
 \le J_D\left(\boldsymbol g\right)-J_D\left(\boldsymbol g_*\right)
 \le\operatorname{gap}_D\left(\boldsymbol g\right).
 \label{eq:cdf-gap-distance}
\end{equation}
In particular, if $\operatorname{gap}_D\left(\boldsymbol g\right)\le\tau$
for some $\tau\ge0$ and $a\le f\le b$ is a density, then
\begin{equation}
 \left\|G_{\boldsymbol g}-F_f\right\|_2
 \le\sqrt\tau+\left\|\widetilde F-F_f\right\|_2+2b/D.
 \label{eq:cdf-gap-risk-pathwise}
\end{equation}
For a measurable feasible output satisfying this optimization error bound almost surely
and deterministic tolerance $\tau$, it follows that
\begin{equation}
 \E_f\left\|G_{\boldsymbol g}-F_f\right\|_2^2
 \le3\left(\tau+\E_f\left\|\widetilde F-F_f\right\|_2^2
                    +4b^2/D^2\right).
 \label{eq:cdf-gap-risk}
\end{equation}
The same conclusions hold for any verified upper bound on
$J_D\left(\boldsymbol g\right)-\min J_D$, in place of the linearization gap.
\end{proposition}

\begin{proof}
The interpolant map is affine and injective, so $J_D$ is a strictly
convex quadratic on the nonempty compact convex set $\mathcal K_D$.
Projection optimality gives
$\left\langle G_{\boldsymbol g_*}-\widetilde F,
G_{\boldsymbol g}-G_{\boldsymbol g_*}\right\rangle\ge0$.
Expanding the squared objective and applying convexity at
$\boldsymbol g$ therefore yields
\[
 \left\|G_{\boldsymbol g}-G_{\boldsymbol g_*}\right\|_2^2
 \le J_D\left(\boldsymbol g\right)-J_D\left(\boldsymbol g_*\right)
 \le\nabla J_D\left(\boldsymbol g\right)\cdot\left(\boldsymbol g-\boldsymbol g_*\right)
 \le\operatorname{gap}_D\left(\boldsymbol g\right).
\]
The true-CDF interpolant $F^D$ is feasible and has
$\left\|F^D-F_f\right\|_2\le b/D$. Nonexpansiveness of projection gives
$\left\|\widehat F_D-F_f\right\|_2\le\left\|\widetilde F-F_f\right\|_2+2b/D$.
Combining this bound with \eqref{eq:cdf-gap-distance} proves
\eqref{eq:cdf-gap-risk-pathwise};
$\left(u+v+w\right)^2\le3\left(u^2+v^2+w^2\right)$ gives \eqref{eq:cdf-gap-risk}.
Only an upper bound on objective suboptimality is needed.
\end{proof}

For the hierarchical step release, $J_D$ is a finite quadratic objective and its linearization gap in \eqref{eq:cdf-optimization-gap} is obtained by linear optimization over $\mathcal K_D$. The gap bound applies to feasible outputs with a verified objective-error bound. Taking $D=2^{\left\lceil\left(\log_2N\right)/2\right\rceil}$ and $\tau\le N^{-1}$ in \eqref{eq:cdf-gap-risk} preserves order $\left(\log N\right)^2/N$. Equations~\eqref{eq:cdf-wasserstein-comparison} and~\eqref{eq:transport-stability} transfer this bound to quantiles and transport maps. This slope-constrained procedure uses $a,b$; Corollary~\ref{cor:quantile} uses the class-parameter-free projection.

\subsection{Finite isotonic projection}
\label{app:finite-isotonic}
Indeed, for $G\in\mathcal D$, its cell averages
$\bar G_j=D\int_{j/D}^{(j+1)/D}G(x)\,dx$ are nondecreasing and
belong to $[0,1]$. On each cell,
\[
 \int_{j/D}^{(j+1)/D}(G-\widehat F_j)^2
 =\int_{j/D}^{(j+1)/D}(G-\bar G_j)^2
       +D^{-1}(\bar G_j-\widehat F_j)^2.
\]
Replacing $G$ by its cell averages therefore preserves feasibility
and cannot increase the objective. This proves the reduction.
The finite objective is strictly convex, so its unique minimizer
is the Euclidean projection onto the displayed closed convex set
and depends continuously on the released values.
\subsection{Quantile and distribution estimation}
\label{app:quantile}

\label{app:cdf-wasserstein-proof}
To see this, let $A\subset\left[0,1\right]^2$ be the region where $t$ lies strictly
between $F_f\left(x\right)$ and $G\left(x\right)$. Up to a set of area zero, it is also the
region where $x$ lies between $q_f\left(t\right)$ and the generalized quantile
$q_P\left(t\right)$. Integrating squared distances in the two directions gives
\[
 \left\|F_f-G\right\|_2^2=2\iint_A|F_f\left(x\right)-t|\,dx\,dt,\qquad
 W_2^2\left(P_f,P\right)=2\iint_A|x-q_f\left(t\right)|\,dx\,dt.
\]
Since $F_f\left(q_f\left(t\right)\right)=t$, the density bounds compare the two integrands.
The argument applies to atomic $P$, including atoms at the endpoints.

Write $q_f=F_f^{-1}$ and let $P_f$ denote the probability measure with
density $f$. For probability measures $P_1,P_2$ on $[0,1]$, monotone coupling gives
$W_2^2\left(P_1,P_2\right)=\int_0^1|F_{P_1}^{-1}\left(t\right)-F_{P_2}^{-1}\left(t\right)|^2\,dt$,
including measures with atoms \citep{PeyreCuturi2018}.

\begin{proof}[Lower bounds in Corollary~\ref{cor:quantile}]
Given any measurable quantile estimate $\widehat q$, clip it pointwise to $\left[0,1\right]$
and project the result in $L^2$ onto $\mathcal C_{1/b,1/a}$, obtaining $G$.
Since $q_f$ lies in this set, these two operations do not increase its
squared $L^2$ error. They are measurable, as are continuous representative
selection and inversion, by the preceding projection argument.
The lower slope $1/b$ of $G$ implies
\[
 |G^{-1}\left(x\right)-F_f\left(x\right)|
 \le b|G\left(F_f\left(x\right)\right)-x|.
\]
Changing variables $x=q_f\left(t\right)$ and using $q_f'\left(t\right)\le1/a$ gives
\[
 \left\|G^{-1}-F_f\right\|_2^2
 \le\frac{b^2}{a}\left\|G-q_f\right\|_2^2
 \le\frac{b^2}{a}\left\|\widehat q-q_f\right\|_2^2.
\]
Thus a quantile estimator induces a CDF estimator by postprocessing the
same public transcript, and the CDF lower bound applies.
For a probability-measure estimate, the second inequality in
\eqref{eq:cdf-wasserstein-comparison} converts its CDF into an estimator
with risk at most $b$ times its squared $W_2$ risk. The CDF lower bound
therefore applies directly, including to atomic outputs.
\end{proof}

\subsection{Transport stability}
\label{app:transport-proof}
Project both CDF estimates by \eqref{eq:cdf-projection} and set $\widehat T=G_Y^{-1}\circ G_X$. For $T=F_g^{-1}\circ F_f$, the lower slope bound gives
\[
 |\widehat T\left(x\right)-T\left(x\right)|
 \le\frac1a|G_X\left(x\right)-G_Y\left(T\left(x\right)\right)|
 \le\frac1a\left(|G_X\left(x\right)-F_f\left(x\right)|
               +|G_Y\left(T\left(x\right)\right)-F_g\left(T\left(x\right)\right)|\right).
\]
Squaring, using $(u+v)^2\le2u^2+2v^2$, and integrating with $T_\#P_f=P_g$ and $f,g\le b$ proves \eqref{eq:transport-stability}. Taking $a=\Gamma^{-1}$ and $b=\Gamma$, \eqref{eq:cdf-projection-risk} and \eqref{eq:cdf-hierarchy-bound} yield the upper bound in \eqref{eq:transport-risk}. The map $\widehat T$ transports the estimated source density to the estimated target density, both satisfying $a\le\widehat f,\widehat g\le b$.

\subsection{Interval anchoring and a positive H\"older prior}

We first verify the interval anchoring identity before constructing the prior.
For any $L^2$ density $f$ on $\left[0,1\right]$, the function $F_f-x$ has equal
zero endpoint values and derivative $f-1$. Subtracting its mean gives
the zero-mean periodic primitive $\mathcal P\left(f-1\right)$. Fubini gives
\[
 \int_0^1\left(F_f\left(x\right)-x\right)\,dx
 =\int_0^1\left(1-u\right)f\left(u\right)\,du-\frac12
 =\frac12-\E_fX,
\]
which proves \eqref{eq:interval-anchor}. Equivalently, integration by
parts yields $\widehat f\left(k\right)=1+2\pi ik\widehat F_f\left(k\right)$ for $k\ne0$;
the boundary term is $F_f\left(1\right)-F_f\left(0\right)=1$.
For two $L^2$ densities $f,g$, the difference $F_f-F_g$ has zero
endpoint values and mean $\E_gX-\E_fX$. Parseval's identity therefore
gives the exact decomposition
\begin{equation}
\label{eq:cdf-periodic-decomposition}
 \left\|F_f-F_g\right\|_2^2
 =\sum_{k\in\mathbb Z\setminus\left\{0\right\}}
       \frac{|\widehat{f-g}\left(k\right)|^2}{\left(2\pi k\right)^2}
       +\left(\E_fX-\E_gX\right)^2.
\end{equation}
The last term records the constant that must accompany a periodic
primitive when estimating an interval CDF.

For $1\le k\le J$, set
\[
 c_k\left(x\right)=\sqrt2\cos\left(2\pi kx\right),\qquad
 s_k\left(x\right)=\sqrt2\sin\left(2\pi kx\right),\qquad
 w_k=\left(2\pi k\right)^{-2},\qquad S_J=\sum_{k=1}^J\sqrt{w_k}.
\]
The $c_k$ have zero integral and are orthonormal in $L^2$; the same
orthonormality holds for the $s_k$. Fix $t>0$ and define
\begin{equation}\label{eq:cdf-prior-model}
 a_k^2=\frac{tS_J}{N_\varepsilon\sqrt{w_k}},\qquad
 f_\theta=1+\sum_{k=1}^J\theta_kc_k,\qquad
 F_\theta\left(x\right)=x+\sum_{k=1}^J\frac{\theta_ks_k\left(x\right)}{2\pi k},
 \quad |\theta_k|\le a_k.
\end{equation}
Equip this rectangle with the product prior
$\Pi=\bigotimes_{k=1}^J\pi_{a_k}$, where $\pi_a$ is the
cosine-squared prior.
Its total prior Fisher information is
\begin{equation}\label{eq:cdf-prior-information}
 \sum_{k=1}^J\frac{\pi^2}{a_k^2}
 =\frac{\pi^2N_\varepsilon}{tS_J}
       \sum_{k=1}^J\sqrt{w_k}
 =\frac{\pi^2N_\varepsilon}{t}.
\end{equation}

For every fixed $s>0$, $\left\|c_k\right\|_{C^s}\le C_sk^s$.
For noninteger $s=m+\alpha$, $0<\alpha<1$, this follows from
the derivative bounds and
$\min\left\{2,2\pi k|x-y|\right\}\le C_\alpha k^\alpha|x-y|^\alpha$.
For integer $s$, the bounds on the derivatives through order $s$
give the assertion directly. Since $S_J\asymp\log J$ for $J\ge2$,
\begin{equation}\label{eq:cdf-holder-cube}
 \sup_{|\theta_k|\le a_k}\left\|f_\theta-1\right\|_{C^s}
 \le C_s\sum_{k=1}^J a_kk^s
 \le C_{s,t}\sqrt{\frac{\log J}{N_\varepsilon}}
                    J^{s+3/2}.
\end{equation}
Choose $J=\left\lfloor N_\varepsilon^\eta\right\rfloor$ with
$0<\eta<\left(2s+3\right)^{-1}$. The last expression tends to zero.
Because $M>\left\|1\right\|_{C^s}$ and $a<1<b$, the entire rectangle in
\eqref{eq:cdf-prior-model} eventually lies in
$\mathcal H_s\left(M,a,b\right)$, with $1/2\le f_\theta\le3/2$.
The parameter-dependent dimension is finite for each experiment, and
the containing H\"older class is fixed.

\subsection{Lower bounds for CDFs and transport maps}
\label{app:cdf-lower}

First consider an arbitrary CDF estimator $H$, allowing any measurable
function output. Infinite Bayes risk already gives the desired lower
bound. Otherwise $H\in L^2$ almost surely under the joint prior and
transcript law, and we may define
\[
 \widetilde\theta_k
   =2\pi k\int_0^1\left[H\left(x\right)-x\right]s_k\left(x\right)\,dx,
 \qquad
 \widehat\theta_k
   =\max\left\{-a_k,\min\left\{a_k,\widetilde\theta_k\right\}\right\}.
\]
Define these coefficients arbitrarily on the null exceptional set.
Bessel's inequality and coordinatewise clipping give
\begin{equation}\label{eq:cdf-coefficient-lower}
 \left\|H-F_\theta\right\|_2^2
 \ge\sum_{k=1}^Jw_k\left(\widetilde\theta_k-\theta_k\right)^2
 \ge\sum_{k=1}^Jw_k\left(\widehat\theta_k-\theta_k\right)^2.
\end{equation}
The subtraction of $x$ uses the primitive in
\eqref{eq:cdf-prior-model}, whose endpoints are $0$ and $1$.

The cosine basis is bounded and orthonormal, and $f_\theta\ge1/2$.
Lemmas~\ref{lem:fisher} and~\ref{lem:van-trees} therefore apply, with
$d_k=\E_\Pi\left[I_{\mathsf Z}\left(\theta\right)\right]_{kk}$:
\[
\E_\Pi\left(\widehat\theta_k-\theta_k\right)^2
\ge\left(d_k+\pi^2/a_k^2\right)^{-1},\qquad
\sum_{k=1}^Jd_k\le2N_\varepsilon.
\]
Weighted Cauchy-Schwarz and \eqref{eq:cdf-prior-information} now give
\begin{equation}\label{eq:cdf-bayes-bound}
 \sum_{k=1}^Jw_k\E_\Pi\left(\widehat\theta_k-\theta_k\right)^2
 \ge\frac{\left(\sum_{k=1}^J\sqrt{w_k}\right)^2}
              {\sum_{k=1}^J\left(d_k+\pi^2/a_k^2\right)}
 \ge\frac{S_J^2}{\left(2+\pi^2/t\right)N_\varepsilon}.
\end{equation}
Combining \eqref{eq:cdf-coefficient-lower} and
\eqref{eq:cdf-bayes-bound}, taking the supremum over the containing
class, and then the infimum over protocols and estimators proves the
lower bound in \eqref{eq:cdf-risk}. Indeed,
$\log J\asymp\log N_\varepsilon\asymp\log N$ and
$N_\varepsilon\asymp N$. This completes the proof of
Theorem~\ref{thm:cdf}.

For transport, restrict the target density to $g=1$ and use the same
source prior. The true map is $T_\theta=F_\theta$, and the potential
\[
 \phi_\theta\left(x\right)=\int_0^x F_\theta\left(u\right)\,du
\]
is twice continuously differentiable with
$\phi_\theta'=F_\theta$ and $\phi_\theta''=f_\theta$.
It is convex, its gradient maps $\left[0,1\right]$ onto $\left[0,1\right]$, and
$\left(\phi_\theta'\right)_\#P_{f_\theta}$ is uniform. By
\eqref{eq:cdf-holder-cube} and the strict inequalities
$\Gamma>1$ and $L_1<1<L_2$, this submodel belongs to
$\mathcal M_s\left(M,\Gamma,L_1,L_2\right)$ for all sufficiently large $N$.

The target users have zero score, so Lemma~\ref{lem:fisher} bounds
the full $2n$-user transcript information by $2N_\varepsilon$.
Equation~\eqref{eq:cdf-bayes-bound} therefore applies to this experiment.
For an arbitrary map
estimator $H$, positivity strengthens \eqref{eq:cdf-coefficient-lower}
to the required comparison
\[
 \left\|H-T_\theta\right\|_{L^2\left(P_{f_\theta}\right)}^2
 \ge\tfrac12\left\|H-F_\theta\right\|_2^2
 \ge\tfrac12\sum_{k=1}^Jw_k\left(\widehat\theta_k-\theta_k\right)^2.
\]
The resulting Bayes bound proves the lower half of
\eqref{eq:transport-risk} and completes
Corollary~\ref{cor:transport}.

\end{document}